\documentclass{article}
\usepackage{amsfonts}
\usepackage{makeidx}
\usepackage{amssymb}
\usepackage{amsmath}
\usepackage{amstext}
\usepackage{graphicx}
\usepackage{latexsym}
\usepackage{amsthm}

\newtheorem{theorem}{Theorem}
\newtheorem{axiom}{Axiom}

\newtheorem{corollary}[theorem]{Corollary}

\newtheorem{definition}[theorem]{Definition}

\newtheorem{lemma}[theorem]{Lemma}

\newtheorem{proposition}[theorem]{Proposition}
\newtheorem{remark}[theorem]{Remark}

\input tcilatex
\begin{document}

\title{The Euclidean Universe}
\author{Vieri Benci \thanks{
Dipartimento di Matematica, Universit\`{a} degli Studi di Pisa, Via F.
Buonarroti 1/c, Pisa, ITALY} \and \bigskip \and \textit{I dedicate this work
to my friend Eugenio who } \and \textit{more than 40 years ago, when we were
two } \and \textit{young guys in New York, he gave me a book } \and \textit{%
by Keisler and introduced me to the beauty } \and \textit{of Nonstandard
Analysis.\bigskip }}
\maketitle

\begin{abstract}
We introduce a mathematical structure called Euclidean Universe. This
structure provides a basic framework for Non-Archimedean Mathematics and
Nonstandard Analysis.

\textbf{Keywords}. Non Archimedean Mathematics, Non Standard Analysis.

\textbf{Mathematics Subject Classification (2020)}: 03H05, 12L15, 03A99. 
\end{abstract}

\tableofcontents

\section{Introduction}

In this paper, we introduce a mathematical structure called Euclidean
Universe. This structure provides a basic framework for Non-Archimedean
Mathematics, namely the Mathematics based on infinite and infinitesimal
numbers.

The Euclidean Universe is defined by three axioms which have been chosen in
such a way to appear absolutely \textit{natural. }The first axiom introduce
the infinite numbers as the \textit{numerosities} of infinite sets in such a
way that the V Euclidean common notion

\begin{center}
\textit{the whole is larger than the part}
\end{center}

\noindent be preserved. In the second axiom we introduce the Euclidean line
with the following peculiarity: if any magnitude can be "represented" by a
point on the Euclidean line, then also the infinite (and consequently the
infinitesimal) magnitudes have this right. Then the Euclidean line must be
larger than the real line. The last axiom is more technical and it is
necessary to make the Euclidean line to include (a copy of) the real numbers.

Actually, this paper can be considered a new introduction to the
Non-Archimedean Mathematics in the spirit of Veronese \cite%
{veronese,veronese2} and Levi-Civita \cite{LC}. Moreover, these axioms are
sufficiently strong to include the basic principles of Nonstandard Analysis
such as the Leibniz Principle and the Saturation Principle. In particular
the Euclidean line incudes, as subfields, infinitely many copies of the
hyperreal numbers of any saturation less or equal to the first inaccessible
number (see Def. \ref{uno}).

Then, the Euclidean Universe includes many results of Non-Archimedean
Mathematics obtained in the last 30 years.

For a recent historical and foundational analysis of the underlying ideas,
we refer to \cite{BF19} and the references therein.

\bigskip

\noindent \textbf{Acknowledgement}: I want to express my gratitude to Prof.
Marco Forti for the many interesting discussions on these topics. Also a big
thank is due to a referee whose remarks made me to improve considerably this
paper.

\section{The three Axioms\label{FA}}

\subsection{The Numerosity axiom\label{NA}}

The first axiom defines the notion of \textbf{numerosity}. The notion of
numerosity was introduced in \cite{BV95,BDN2003,BDNF1} as a generalization
of finite cardinality that also applies to infinite sets. The main feature
of numerosities is that they preserve the spirit of the ancient Euclidean
principle that \textit{the whole is larger than the part}; indeed, the
numerosity of a proper subset is strictly smaller than the numerosity of the
whole set.

In principle it would be desirable to define the numerosity for the class of
all sets; however, in order to develop the theory, it is convenient to work
in a "universe" which is a set closed with respect to the main set
operations provided that it is very large. In order to do this we recall a
well known notion in set theory:

\begin{definition}
\label{uno}A cardinal number $\chi $ is \textbf{inaccessible} if it is not a
sum of fewer than $\chi $ cardinals that are less than $\chi $ and $\zeta
<\chi $ implies $2^{\zeta }<\chi $.\ $\chi $ is \textbf{strongly
inaccessible\ }if it is inaccessible and uncountable.
\end{definition}

The first inaccessible cardinal number is $\aleph _{0}.$ The first strongly
inaccessible cardinal number will be denoted by $\kappa .$ The existence of
sets of strongly inaccessible cardinality is established by the Axiom of
Inaccessibility which is independent from ZFC. We will assume this axiom and
in particular we will assume that there exists a set of atoms\footnote{%
In set theory, an atom $a$ is any entity that is not a set, namely $a$ is an
atom if and only if%
\begin{equation*}
\forall x,\ x\notin a
\end{equation*}%
} $\mathbb{A}$ having cardinality $\kappa .$ Then we can define a "universe" 
$\Lambda $ defined as follows:%
\begin{equation}
\Lambda =\left\{ X\in V(\mathbb{A})\ |\ \left\vert X\right\vert <\kappa
\right\}  \label{l1}
\end{equation}%
where for any set $A,$ $V(A)$ denotes the superstructure over $A$ namely%
\begin{equation*}
V(A)=\bigcup_{n\in \mathbb{N}}V_{n}(A)
\end{equation*}%
with%
\begin{equation*}
V_{0}(A)=A
\end{equation*}%
and, for every $n\in \mathbb{N}$, 
\begin{equation}
V_{n+1}(A)=V_{n}(A)\cup \mathfrak{\wp }\left( V_{n}(A)\right) .  \label{Vn}
\end{equation}%
We will refer to $\Lambda $ as to the \textbf{accessible universe} since its
sets have strongly accessible cardinality (and finite rank\footnote{%
The rank of \ $\varnothing \ $is $0$.\ The rank of a set $E\neq 0$ is the
least ordinal number greater than the rank of any element of the set.}). $%
\Lambda $ can be split as follows%
\begin{equation*}
\Lambda =\Lambda _{S}\cup \mathbb{A}
\end{equation*}%
where $\Lambda _{S}$ is a family of sets. Notice that $\mathbb{A}\notin
\Lambda $ since $\left\vert \mathbb{A}\right\vert =\kappa .$

Now we are ready to state our first axiom:

\begin{axiom}
\label{num}\textbf{(Numerosity axiom) }The numerosity is a surjective map%
\begin{equation*}
\mathfrak{num}:\Lambda _{S}\rightarrow \mathfrak{N},\ \ \ \mathfrak{N}%
\subset \mathbb{A}
\end{equation*}%
which satisfies the following properties: if $a,b\in \Lambda $ and $A$,$B$,$%
A^{\prime }$,$B^{\prime }\in \Lambda _{S},$

\begin{enumerate}
\item \label{n0}$\mathfrak{num}\left( \left\{ a\right\} \right) =\mathfrak{%
num}\left( \left\{ b\right\} \right) $

\item \label{n1}if $A\subset B\ $strictly, then%
\begin{equation*}
\mathfrak{num}\left( A\right) <\mathfrak{num}\left( B\right)
\end{equation*}

\item \label{n3}if $A\cap B=\varnothing ,$\ $\mathfrak{num}\left( A\right) =%
\mathfrak{num}\left( A^{\prime }\right) ,\ \mathfrak{num}\left( B\right) =%
\mathfrak{num}\left( B^{\prime }\right) ,$ then%
\begin{equation*}
\mathfrak{num}\left( A\cup B\right) =\mathfrak{num}\left( A^{\prime }\cup
B^{\prime }\right)
\end{equation*}

\item \label{n4}if $\mathfrak{num}\left( A\right) =\mathfrak{num}\left(
A^{\prime }\right) ,\ \mathfrak{num}\left( B\right) =\mathfrak{num}\left(
B^{\prime }\right) ,$ then%
\begin{equation*}
\mathfrak{num}\left( A\times B\right) =\mathfrak{num}\left( A^{\prime
}\times B^{\prime }\right)
\end{equation*}

\item \label{n5}if $A\in \Lambda _{S}$ and $b\in \Lambda ,$ then $\mathfrak{%
num}\left( A\times \left\{ b\right\} \right) =\mathfrak{num}\left( A\right)
. $
\end{enumerate}
\end{axiom}

\bigskip

If $F$ and $F^{\prime }$ are finite sets of the same cardinality, by \ref%
{num}.\ref{n0} and \ref{num}.\ref{n3}, it follows that $\mathfrak{num}\left(
F\right) =\mathfrak{num}\left( F^{\prime }\right) ;$ then, by \ref{num}.\ref%
{n5}, the numerosities of finite sets can be identified with the natural
numbers $\mathbb{N}$. By \ref{num}.\ref{n3}, $\mathfrak{N}$ can be equipped
with an "addition" by setting%
\begin{equation*}
\sigma +\tau =\mathfrak{num}\left( A\cup B\right)
\end{equation*}%
where $\sigma =\mathfrak{num}\left( A\right) ,$ $\tau =\mathfrak{num}\left(
B\right) $ and $A\cap B=\varnothing $; similarly, by \ref{num}.\ref{n4}, we
can define a "multiplication":%
\begin{equation*}
\sigma \cdot \tau =\mathfrak{num}\left( A\times B\right)
\end{equation*}%
Clearly $0=\mathfrak{num}\left( \varnothing \right) $ is the neutral element
with respect to the addition and, by \ref{num}.\ref{n5}, $\mathfrak{num}%
\left( \left\{ b\right\} \right) =1$ is the neutral element with respect to
the multiplication for any $b\in \Lambda .$

\subsection{The Euclidean Field Axiom\label{EFA}}

Our second axioms identifies the Euclidean (straight) line with a field $%
\mathbb{E}_{\kappa }\subseteq \mathbb{A}$. Usually the Euclidean line is
identified with the real line, however we think that this point of view is
too restrictive. In fact, the main intuitive peculiarities of the Euclidean
line are the following:

\begin{itemize}
\item two oriented segments of the Euclidean line can be added or subtracted;

\item if we choose a unitary segment of the Euclidean line, two segments can
be multiplied or divided;

\item once we have chosen two distinguish points $O$ and $U$ on the line,
every magnitude can be posed in a biunivocal correspondence with a point
(provided that the unitary magnitude has been defined).
\end{itemize}

Then, if we take $0=O$ and $1=U,$ the Euclidean line gets the stucture of
ordered field and its points can be identified with numbers. Since the
numerosities can be considered magnitudes, the Euclidean line should be
richer than the real line. We can formalize these intuitive remarks by the
following axiom:

\begin{axiom}
\label{1}(\textbf{Euclidean} \textbf{Field Axiom}) There is an ordered field 
$\mathbb{E}_{\kappa }\subset \mathbb{A}$ such that

\begin{itemize}
\item $\mathbb{E}_{\kappa }$ contains the set numerosities $\mathfrak{N}$
and the field operations $+,\cdot $ coincide with the numerosity operations;

\item for every $\xi \in \mathbb{E}_{\kappa },\ $there exists $E\in \Lambda
, $ such that%
\begin{equation}
\left\vert \xi \right\vert <\mathfrak{num}\left( E\right) .  \label{bbb}
\end{equation}
\end{itemize}
\end{axiom}

We will refer to $\mathbb{E}_{\kappa }$ as the \textbf{Euclidean line} or
the \textbf{Euclidean} \textbf{field} and its elements will be called 
\textbf{Euclidean numbers}. $\mathbb{E}_{\kappa }$ contains infinite
numbers, then it is a non-Archimedean field. Now let us recall some basic
definitions relative to non-Archimedean fields. Since $\mathbb{N\subset E}%
_{\kappa }$,\ the following definition makes sense:

\begin{definition}
\label{in}Let $\xi \in \mathbb{E}_{\kappa }$. We say that:

\begin{itemize}
\item $\xi $ is \textbf{infinitesimal} if $\forall n\in \mathbb{N}$, $|\xi |<%
\frac{1}{n}$;

\item $\xi $ is \textbf{finite }(or bounded) if $\exists n\in \mathbb{N}$
such as $|\xi |<n$;

\item $\xi $ is \textbf{infinite} (or unbounded) if $\forall n\in \mathbb{N}$%
, $|\xi |>n.$
\end{itemize}
\end{definition}

The following proposition establishes some relations among Euclidean numbers:

\begin{proposition}
\label{t235}We have the following relations:

\begin{itemize}
\item (i) If $\varepsilon $ e ${\delta }$ are infinitesimal, also $%
\varepsilon +{\delta }$ ,$\mathbf{\;}\varepsilon -{\delta }$ and $%
\varepsilon \cdot {\delta }$ are infinitesimal.

\item (ii) If $\xi $ e $\sigma $ are bounded also $\xi +\sigma $ ,$\mathbf{\;%
}\xi -\sigma $ e $\xi \cdot \sigma $ are bounded.

\item (iii) If $\theta $ e $\tau $ are infinite, also $\theta \cdot \tau $
is infinite; moreover if $\theta $ and $\tau $ are postive infinite (or
negative), also $\omega +\tau $ is postive infinite (or negative).

\item (iv) If $\varepsilon $ is infinitesimal and $\xi $ is bouded $%
\varepsilon \cdot \xi $ is infinitesimal; moreover if $\varepsilon \neq 0$
and $\xi $ is not infinitesimal, $\xi /\varepsilon $ is infinite.
\end{itemize}
\end{proposition}

\textbf{Proof. }The proof of this proposition is easy and it is left to the
reder.

$\square $\bigskip

\begin{definition}
\label{dic} We say that two numbers $\xi ,\zeta \in \mathbb{E}${$_{\kappa }$}
are \textbf{infinitely close} if $\xi -\zeta $ is infinitesimal. In this
case, we write $\xi \sim \zeta $.
\end{definition}

Clearly, the relation "$\sim $" of infinite closeness is an equivalence
relation. Then the following definition comes naturally

\begin{definition}
If $\xi \in \mathbb{E}_{\kappa }$, the monad of $\xi $ is the set of all
numbers that are infinitely close to it:%
\begin{equation*}
\mathfrak{m}\mathfrak{o}\mathfrak{n}(\xi )=\{\zeta \in \mathbb{E}_{\kappa }%
\mathbb{\ }|\ \xi \sim \zeta \},
\end{equation*}%
The galaxy of $\xi $ is the set of all numbers that are finitely close to it:%
\begin{equation*}
\mathfrak{gal}(\xi )=\{\zeta \in \mathbb{E}_{\kappa }\mathbb{\ }|\ \xi
-\zeta \text{\ is a finite number}\}.
\end{equation*}
\end{definition}

\subsection{The Center Axiom}

The notion of monad allows to state our last axiom:

\begin{axiom}
\label{3}Every monad $\mu $ has a distinguished point called \textbf{center}
of $\mu $ and denoted by $Ctr(\mu )$; the set $\mathfrak{C}$ of all the
centers is an additive subgroup of $\mathbb{E}_{\kappa }$ containing $%
\mathfrak{N.}$
\end{axiom}

For every $\xi \in \mathbb{E}_{\kappa }$, we will use the notation%
\begin{equation*}
ctr\left( \xi \right) =Ctr(\mathfrak{mon}\left( \xi \right) );
\end{equation*}%
so, every number $\xi \in \mathbb{E}_{\kappa }$ can be decomposed as follows:%
\begin{equation}
\xi =x+\varepsilon  \label{split}
\end{equation}%
where $x=ctr\left( \xi \right) \ $and $\varepsilon \sim 0$. Then $\mathbb{E}%
_{\kappa }$ can be decomposed as follows:%
\begin{equation*}
\mathbb{E}_{\kappa }=\mathfrak{C}\times \mathfrak{mon}\left( 0\right)
\end{equation*}%
and%
\begin{equation*}
ctr:\mathbb{E}_{\kappa }\rightarrow \mathfrak{C}
\end{equation*}%
is a projection.

\section{Structure of the Euclidean numbers}

In this section we will examine some peculiarities of $\mathbb{E}_{\kappa }.$
In particular, we will see that $\mathbb{E}_{\kappa }$ contains the ordinal
numbers of accessible cardinalities and the real numbers. Moreover we will
introduce the notion $\Lambda $-limit which, in this context, is a very
basic tool.

\subsection{The ordinal numerosities\label{ON}}

In this subsection we introduce a set $\mathbf{Ord}\subset \mathfrak{N}$
that is isomorphic (in a sense specified below) to the initial segment of
length $\kappa $ of ordinal numbers.

\begin{definition}
The set $\mathbf{Ord}\subset \mathfrak{N}$ of ordinal numerosities is
defined as follows: $\tau \in \mathbf{Ord}$ if and only if%
\begin{equation*}
\tau =\mathfrak{num}\left( \Omega _{\tau }\right) ,
\end{equation*}%
where 
\begin{equation*}
\Omega _{\tau }=\left\{ x\in \mathbf{Ord\ |\ }x<\tau \right\} .
\end{equation*}
\end{definition}

It is easy to see by transfinite induction that this is a good recursive
definition. In fact, it is immediate to check that

\begin{itemize}
\item $0\in \mathbf{Ord}$;

\item if $\tau \in \mathbf{Ord,}$ then $\tau +1=\mathfrak{num}\left( \Omega
_{\tau }\cup \left\{ \tau \right\} \right) \in \mathbf{Ord}$ (and hence $%
\mathbb{N}\subset \mathbf{Ord}$).

\item if $\tau _{k}=\mathfrak{num}\left( \Omega _{k}\right) ,$ $k\in K,$ ($%
\left\vert K\right\vert <\kappa $) are ordinal numerosities, then 
\begin{equation*}
\tau :=\mathfrak{num}\left( \dbigcup\limits_{k\in K}\Omega _{k}\right) \in 
\mathbf{Ord.}
\end{equation*}
\end{itemize}

In particular, $\omega =\mathfrak{num}\left( \mathbb{N}\right) $ is an
ordinal. This construction of the ordinal numbers is similar to the
construction of Von Neumann. However, whilst a Von Neumann ordinal $\tau $
is the set of all the Von Neumann ordinals contained in $\tau $, in our
construction an ordinal $\tau $ is the numerosity of the set of ordinals
smaller than $\tau $. Hence, here, an ordinal number, as any other
numerosity, is an atom in $\mathbb{E}_{\kappa }$.

It is easy to see that $\mathbf{Ord}$ is well ordered and hence it is
isomorphic to the initial segment of length $\kappa $ of the full class of
the ordinals. Obviously, not all numerosities are ordinals: for example, $%
\omega -1=\mathfrak{num}\left( \mathbb{N}^{+}\right) =\mathfrak{num}\left( 
\mathbb{N}\backslash \left\{ 0\right\} \right) $ is not an ordinal.

The most remarkable thing in this theory is that the numerosity operations $%
+ $ and $\cdot ,$ correspond to the \textit{natural} (or Hessenberg)
operations between ordinals. We refer to \cite{BLB21} for an in-depth
analysis of this topic.

\begin{remark}
Notice that the existence of $\mathbf{Ord}$ depends only on Axiom 1.
\end{remark}

\subsection{The $\Lambda $-limit theorem\label{LLL}}

We set 
\begin{equation*}
\mathfrak{L}=\left\{ \lambda \in \Lambda \ |\ \lambda \ \text{is a finite set%
}\right\}
\end{equation*}%
and%
\begin{equation*}
\mathfrak{F}\left( \mathfrak{L},\mathbb{E}_{\kappa }\right) =\left\{ \varphi
\in \mathbb{E}_{\kappa }^{\mathfrak{L}}\ |\ \exists A\in \Lambda ,\forall
\lambda \in \mathfrak{L},\ \varphi (\lambda )=\varphi (\lambda \cap
A)\right\} .
\end{equation*}%
Since $\left( \mathfrak{L},\subseteq \right) $ is a directed sets, the
elements of $\mathfrak{F}\left( \mathfrak{L},\mathbb{E}_{\kappa }\right) $
are nets. The set $\mathfrak{F}\left( \mathfrak{L},\mathbb{E}_{\kappa
}\right) $ is a partially ordered commutative algebra over $\mathbb{E}%
_{\kappa }$ with the operations%
\begin{eqnarray*}
\left( \varphi +\psi \right) (\lambda ) &=&\varphi (\lambda )+\psi (\lambda )
\\
\left( \varphi \cdot \psi \right) (\lambda ) &=&\varphi (\lambda )\cdot \psi
(\lambda )
\end{eqnarray*}

\begin{theorem}
\label{LL}($\mathbf{\Lambda }$\textbf{-limit theorem}) There is a unique
ring homomorphism%
\begin{equation*}
J:\mathfrak{F}\left( \mathfrak{L},\mathbb{E}_{\kappa }\right) \rightarrow 
\mathbb{E}_{\kappa }
\end{equation*}%
such that,%
\begin{equation*}
\forall A\in \Lambda ,\ J(\psi _{A})=\mathfrak{num}\left( A\right)
\end{equation*}%
where%
\begin{equation}
\psi _{A}\left( \lambda \right) =\left\vert A\cap \lambda \right\vert .
\label{bella}
\end{equation}
\end{theorem}

\textbf{Proof}: Let $\mathfrak{F}_{q}\left( \mathfrak{L},\mathbb{E}_{\kappa
}\right) $ be the $\mathbb{E}_{\kappa }$-subalgebra of $\mathfrak{F}\left( 
\mathfrak{L},\mathbb{E}_{\kappa }\right) $ "generated" by $\left\{ \psi
_{A}\ |\ A\in \Lambda \right\} \ $namely the subset of the elements $\varphi 
$ of $\mathfrak{F}\left( \mathfrak{L},\mathbb{E}_{\kappa }\right) $ which
can be written as follows:%
\begin{equation*}
\varphi (\lambda )=\frac{\sum_{A\in \mathcal{A}}a_{A}\psi _{A}(\lambda )}{%
\sum_{B\in \mathcal{B}}b_{A}\psi _{B}(\lambda )}
\end{equation*}%
where $\mathcal{A},\mathcal{B}$ are finite subsets of $\Lambda $, $%
a_{A},b_{A}\in \mathbb{E}_{\kappa }$, $\psi _{A},\psi _{B}$ are defined by (%
\ref{bella}) and $\ \forall \lambda \in \mathfrak{L,\ }\sum_{B\in \mathcal{B}%
}b_{A}\psi _{B}(\lambda )\neq 0$. We define a field homomorphism 
\begin{equation*}
J_{q}:\mathfrak{F}_{q}\left( \mathfrak{L},\mathbb{E}_{\kappa }\right)
\rightarrow \mathbb{E}_{\kappa }
\end{equation*}%
as follows:$\ $%
\begin{equation*}
J_{q}(\varphi )=\frac{\sum_{A\in \mathcal{A}}a_{A}\cdot \mathfrak{num}\left(
A\right) }{\sum_{B\in \mathcal{B}}b_{A}\cdot \mathfrak{num}\left( B\right) }
\end{equation*}%
Since $\func{Im}\left( J_{q}\right) =\mathbb{E}_{\kappa }$ is a field, $\ker
\left( J_{q}\right) $ is a maximal ideal in $\mathfrak{F}_{q}\left( 
\mathfrak{L},\mathbb{E}_{\kappa }\right) ;$ hence, the set 
\begin{equation*}
\mathcal{U}_{0}=\left\{ Q\in \mathfrak{L}\ |\ \exists \psi \in \ker \left(
J_{q}\right) ,\ Q=\psi ^{-1}(0)\right\}
\end{equation*}%
is a filter over $\mathfrak{L}$. We denote by $\mathcal{U}$ an ultrafilter
such that $\mathcal{U}_{0}\subseteq \mathcal{U}$. Then also 
\begin{equation*}
\mathfrak{I}:=\left\{ \psi \in \mathfrak{F}\left( \mathfrak{L},\mathbb{E}%
_{\kappa }\right) \ |\ \exists Q\in \mathcal{U},\ \forall \lambda \in Q,\
\psi (\lambda )=0\right\}
\end{equation*}%
is a maximal ideal in $\mathfrak{F}\left( \mathfrak{L},\mathbb{E}_{\kappa
}\right) $ and hence 
\begin{equation*}
\mathbb{F}:=\mathfrak{F}\left( \mathfrak{L},\mathbb{E}_{\kappa }\right) /%
\mathfrak{I}
\end{equation*}%
is a field. We will see that $\mathbb{F}$ is isomorphic to $\mathbb{E}%
_{\kappa }.$ We denote by $\left[ \varphi \right] $ a generic element of $%
\mathbb{F}$ and we claim that%
\begin{equation}
\forall \left[ \varphi \right] \in \mathbb{F},\ \exists \bar{\mu}\in 
\mathfrak{L},\ \exists Q\in \mathcal{U},\ \forall \lambda \in Q,\ \varphi
\left( \lambda \right) =\varphi \left( \bar{\mu}\right) ,  \label{traviata}
\end{equation}%
namely 
\begin{equation}
\left[ \varphi \right] =\left[ C_{\varphi \left( \bar{\mu}\right) }\right]
\label{tosca}
\end{equation}%
where $\lambda \mapsto C_{\xi }\left( \lambda \right) $ denotes the net
identically equal to $\xi .$ In order to prove (\ref{tosca}) we set%
\begin{equation}
R^{-}:=\left\{ \mu \in \mathfrak{L}\ |\ \left[ C_{\varphi (\mu )}\right] <%
\left[ \varphi \right] \right\}  \label{piu}
\end{equation}%
\begin{equation}
R^{0}:=\left\{ \mu \in \mathfrak{L}\ |\ \left[ C_{\varphi (\mu )}\right] =%
\left[ \varphi \right] \right\}  \label{zero}
\end{equation}%
\begin{equation}
R^{+}:=\left\{ \mu \in \mathfrak{L}\ |\ \left[ C_{\varphi (\mu )}\right] >%
\left[ \varphi \right] \right\}  \label{meno}
\end{equation}%
By (\ref{piu}), if $R^{-}\neq \varnothing ,\ \forall \mu \in R^{-},\ \exists
Q_{\mu }^{-}\in \mathcal{U}$ such that 
\begin{equation*}
\forall \lambda \in Q_{\mu }^{-},\ C_{\varphi (\mu )}\left( \lambda \right)
<\varphi \left( \lambda \right)
\end{equation*}%
then,%
\begin{equation*}
\mu \in R^{-}\cap Q_{\mu }^{-}\Rightarrow \ C_{\varphi (\mu )}\left( \mu
\right) <\varphi \left( \mu \right)
\end{equation*}%
and since, by definition $C_{\varphi (\mu )}\left( \mu \right) =\varphi (\mu
),$ it follows that 
\begin{equation*}
\forall \mu \in R^{-},\ R^{-}\cap Q_{\mu }^{-}=\varnothing
\end{equation*}%
and hence,%
\begin{equation*}
R^{-}\notin \mathcal{U}
\end{equation*}%
By (\ref{meno}), arguing in the same way, we have that 
\begin{equation*}
R^{+}\notin \mathcal{U}
\end{equation*}%
Since $\left( R^{-}\cup R^{+}\right) \cup R^{0}=\mathfrak{L}$, it follows
that $R^{0}\in \mathcal{U}$ and hence $R^{0}\neq \varnothing $. Now, if you
take $\bar{\mu}$ in $R^{0}$, there is $Q^{0}\in \mathcal{U}$ such that 
\begin{equation*}
\forall \lambda \in Q^{0},\ \varphi \left( \lambda \right) =C_{\varphi (\bar{%
\mu})}\left( \lambda \right) =\varphi (\bar{\mu}).
\end{equation*}%
namely (\ref{tosca}) is satisfied. Now, we can extend $J_{q}$ to $\mathfrak{F%
}\left( \mathfrak{L},\mathbb{E}_{\kappa }\right) ;$ given $\varphi \in 
\mathfrak{F}\left( \mathfrak{L},\mathbb{E}_{\kappa }\right) $, using (\ref%
{traviata}) we set%
\begin{equation*}
J\left( \varphi \right) =\left[ C_{\varphi \left( \bar{\mu}\right) }\right] .
\end{equation*}%
So every function $\varphi $ in $\mathfrak{F}\left( \mathfrak{L},\mathbb{E}%
_{\kappa }\right) $ is eventually constant in the sense that%
\begin{equation*}
\exists \xi \in \mathbb{E}_{\kappa },\exists Q\in \mathcal{U},\forall
\lambda \in Q,\ \varphi \left( \lambda \right) =\xi .
\end{equation*}%
Then $\mathbb{F}$ is isomorphic to $\mathbb{E}_{\kappa }.$

It remains to prove the uniqueness of $J.$ Let us assume that $J_{1}$ and $%
J_{2}$ extend $J_{q}$ to all $\mathfrak{F}\left( \mathfrak{L},\mathbb{E}%
_{\kappa }\right) .$ We have to prove that for every $\varphi \in \mathfrak{F%
}\left( \mathfrak{L},\mathbb{E}_{\kappa }\right) $%
\begin{equation*}
J_{1}\left( \varphi \right) =J_{2}\left( \varphi \right)
\end{equation*}

We set%
\begin{equation*}
c_{1}=J_{1}\left( \varphi \right) ;\ c_{2}=J_{2}\left( \varphi \right)
\end{equation*}%
\begin{eqnarray*}
A_{\varphi } &=&\left\{ J_{q}\left( \psi \right) \ |\ \psi \in \mathfrak{F}%
_{q}\left( \mathfrak{L},\mathbb{E}_{\kappa }\right) ,\ \exists Q\in \mathcal{%
U}_{0},\ \forall \lambda \in Q,\ \psi (\lambda ),\ \psi \leq \varphi
\right\} ;\  \\
B_{\varphi } &=&\left\{ J_{q}\left( \psi \right) \ |\psi \in \mathfrak{F}%
_{q}\left( \mathfrak{L},\mathbb{E}_{\kappa }\right) \ \ \exists Q\in 
\mathcal{U}_{0},\ \forall \lambda \in Q,\ \psi (\lambda ),\ \psi >\varphi
\right\}
\end{eqnarray*}

Clearly $\ \forall a\in A_{\varphi },\ \forall b\in B_{\varphi },\ $ 
\begin{equation}
a\leq c_{1}\leq b\ \text{\ \ and\ \ }\ a\leq c_{2}\leq b.  \label{ccc}
\end{equation}%
and hence, assuming that $c_{1}\leq c_{2}$ 
\begin{equation*}
\forall a\in A_{\varphi },\ \forall b\in B_{\varphi },\ 0\leq
c_{2}-c_{1}\leq b-a
\end{equation*}%
Since $A_{\varphi }\cup B_{\varphi }=\mathfrak{F}_{q}\left( \mathfrak{L},%
\mathbb{E}_{\kappa }\right) $ contains $\frac{1}{\mathfrak{num}\left(
E\right) }$ for any set $E\in \Lambda _{S}\backslash \left\{ 0\right\} $, we
have that 
\begin{equation*}
\ 0\leq c_{2}-c_{1}\leq \frac{1}{\mathfrak{num}\left( E\right) }
\end{equation*}%
and so, by (\ref{bbb}), $c_{2}=c_{1}.$

$\square $\bigskip

\begin{definition}
\label{lim+}The number $J\left( \varphi \right) $ is called $\Lambda $-limit
of the net $\varphi $ and will be denoted by%
\begin{equation*}
J\left( \varphi \right) =\lim_{\lambda \uparrow \Lambda }\varphi (\lambda ).
\end{equation*}
\end{definition}

The reason of this name and notation is that the operation%
\begin{equation*}
\varphi \mapsto \lim_{\lambda \uparrow \Lambda }\varphi (\lambda )
\end{equation*}%
satisfies some of the properties of the usual limit over a net:

\begin{itemize}
\item If eventually $\varphi (\lambda )\geq \psi (\lambda )$, then%
\begin{equation*}
\lim_{\lambda \uparrow \Lambda }\ \varphi (\lambda )\geq \lim_{\lambda
\uparrow \Lambda }~\psi (\lambda ).
\end{equation*}

\item If $\forall q\in \mathbb{Q}$, $C_{q}\left( \lambda \right) =q$, then%
\begin{equation*}
\lim_{\lambda \uparrow \Lambda }\ C_{q}(\lambda )=q.
\end{equation*}

\item For all $\varphi ,\psi \in \mathfrak{F}\left( \mathfrak{L},\mathbb{E}%
_{\kappa }\right) $ 
\begin{eqnarray*}
\lim_{\lambda \uparrow \Lambda }\varphi (\lambda )+\lim_{\lambda \uparrow
\Lambda }\psi (\lambda ) &=&\lim_{\lambda \uparrow \Lambda }\left( \varphi
(\lambda )+\psi (\lambda )\right) , \\
\lim_{\lambda \uparrow \Lambda }\varphi (\lambda )\cdot \lim_{\lambda
\uparrow \Lambda }\psi (\lambda ) &=&\lim_{\lambda \uparrow \Lambda }\left(
\varphi (\lambda )\cdot \psi (\lambda )\right) .
\end{eqnarray*}
\end{itemize}

In this framework, $\Lambda $ can be regarded as the "point at infinity" of $%
\mathfrak{L}$. The $\Lambda $-limit is not a limit in a topological sense,
in fact there are also strong differences with a topological limit; we list
some of them:

\begin{itemize}
\item Every net $\varphi \in \mathfrak{F}\left( \mathfrak{L},\mathbb{E}%
_{\kappa }\right) $ has a limit $L\in \mathbb{E}_{\kappa }$.

\item If $\forall \lambda \in \mathfrak{L}$, $\varphi (\lambda )\neq 0$,
then $\lim_{\lambda \uparrow \Lambda }\varphi (\lambda )\neq 0;$ in fact,%
\begin{equation*}
\lim_{\lambda \uparrow \Lambda }\varphi (\lambda )\cdot \lim_{\lambda
\uparrow \Lambda }\frac{1}{\varphi (\lambda )}=1
\end{equation*}%
and hence $\lim_{\lambda \uparrow \Lambda }\varphi (\lambda )\neq 0.$

\item If, $\xi \in \mathbb{E}_{\kappa }\backslash \mathbb{R}$,%
\begin{equation*}
\lim_{\lambda \uparrow \Lambda }C_{\xi }\left( \lambda \right) \neq \xi
\end{equation*}
\end{itemize}

For example, take 
\begin{equation*}
\omega :=\lim_{\lambda \uparrow \Lambda }\left\vert \lambda \cap \mathbb{N}%
\right\vert ;
\end{equation*}%
then $\forall \lambda \in \mathfrak{L,\ }\left\vert \lambda \cap \mathbb{N}%
\right\vert <\omega ,$ so 
\begin{equation*}
0>\lim_{\lambda \uparrow \Lambda }\left( \left\vert \lambda \cap \mathbb{N}%
\right\vert -\omega \right) =\lim_{\lambda \uparrow \Lambda }\left\vert
\lambda \cap \mathbb{N}\right\vert -\lim_{\lambda \uparrow \Lambda }\omega
=\omega -\lim_{\lambda \uparrow \Lambda }\omega .
\end{equation*}%
and hence%
\begin{equation*}
\omega <\lim_{\lambda \uparrow \Lambda }\omega .
\end{equation*}

The last statement suggests the following notation: for any $\xi \in \mathbb{%
E}_{\kappa }$, we set 
\begin{equation}
\xi ^{\ast }=\lim_{\lambda \uparrow \Lambda }\xi .  \label{star-}
\end{equation}

\subsection{The real numbers}

We remark that the notion and the first properties of the $\Lambda $-limit
do not depend on Axiom \ref{3}. In this section we will see that also Axiom %
\ref{3} is very relevant.

\begin{definition}
An Euclidean number is called \textbf{standard} if it is finite and it is
the center of a monad. The set of standard points will be denoted by $%
\mathbb{R}$, namely%
\begin{equation*}
\mathbb{R}:=\mathfrak{C}\cap \mathfrak{gal}(0).
\end{equation*}%
If $\xi $ is a finite number, then $ctr(\xi )$ is called \textbf{standard
part }of $x$ and will be denoted also by $st\left( \xi \right) .$
\end{definition}

First let us examine some (obvious) properties of the function $st(\cdot ).$

\begin{proposition}
\label{PS}Let $\xi $ and $\zeta $ be finite numbers, then

\begin{enumerate}
\item \label{blo}$\xi \in \mathbb{R}\Leftrightarrow st\left( \xi \right)
=\xi ;$

\item \label{d}$\xi \leq \zeta \Rightarrow st\left( \xi \right) \leq
st\left( \zeta \right) ;$

\item $st\left( \xi +\zeta \right) =st\left( \xi \right) +st\left( \zeta
\right) ;$

\item \label{4}$st\left( \xi \cdot \zeta \right) =st\left( \xi \right) \cdot
st\left( \zeta \right) ;$

\item \label{55}if $st\left( \zeta \right) \neq 0,$ then $st\left( \frac{\xi 
}{\zeta }\right) =\frac{st\left( \xi \right) }{st\left( \zeta \right) }.$
\end{enumerate}
\end{proposition}

\textbf{Proof}: The first four statements trivially descend from (\ref{split}%
) and Prop. \ref{t235}. In order to prove \ref{PS}.\ref{4}, we put%
\begin{eqnarray*}
\xi &=&r+\varepsilon \\
\zeta &=&s+\theta
\end{eqnarray*}%
where $r,s\in \mathbb{R}$, $\varepsilon \sim \theta \sim 0.$ Then, 
\begin{equation*}
st\left( \xi \cdot \zeta \right) =st\left[ \left( r+\varepsilon \right)
\left( s+\theta \right) \right] =st\left[ rs+\left( \varepsilon s+\theta
r+\varepsilon \theta \right) \right]
\end{equation*}%
Since $\varepsilon s+\theta r+\varepsilon \theta \sim 0$, we have that $%
st\left( \xi \cdot \zeta \right) =rs=st\left( \xi \right) \cdot st\left(
\zeta \right) .$ Let us prove \ref{PS}.\ref{55}; by \ref{PS}.\ref{4} we have
that 
\begin{equation*}
st\left( \zeta \right) \cdot st\left( \frac{\xi }{\zeta }\right) =st\left(
\zeta \cdot \frac{\xi }{\zeta }\right) =st\left( \xi \right) ;
\end{equation*}%
hence%
\begin{equation*}
st\left( \frac{\xi }{\zeta }\right) =\frac{st\left( \xi \right) }{st\left(
\zeta \right) }.
\end{equation*}

$\square $\bigskip

\begin{theorem}
\label{sat}The set of standard numbers $\mathbb{R}$ is isomorphic to the set
of real numbers.
\end{theorem}

\textbf{Proof}. We will prove that every Cauchy sequence of rationals is
convergent to some $L\in \mathbb{R}$ with respect to the metric topology.
Let $x_{n}$ be a Cauchy sequence in $\mathbb{Q}$. We set%
\begin{equation*}
\varphi (\lambda ):=x_{\left\vert \mathbb{N}\cap \lambda \right\vert }
\end{equation*}%
and 
\begin{equation*}
L=st\left( \lim_{\lambda \uparrow \Lambda }\ \varphi \left( \lambda \right)
\right) .
\end{equation*}%
Then, by Prop. \ref{PS}.\ref{blo}, $L\in \mathbb{R}$. We have to prove that $%
L$ is the Cauchy limit of $x_{n}$. We choose a number $\varepsilon \in 
\mathbb{Q}^{+}$; then, there exists $n_{0}$ such that $\forall n,m\geq
n_{0}, $ 
\begin{equation*}
\left\vert x_{n}-x_{m}\right\vert <\varepsilon
\end{equation*}%
Now take $\lambda _{0}\in \mathfrak{L}$ such that $\left\vert \mathbb{N}\cap
\lambda _{0}\right\vert \geq n_{0};\ $thus, $\forall \lambda \supset \lambda
_{0},\ $we have that $\left\vert \mathbb{N}\cap \lambda \right\vert \geq
n_{0}.$ Then%
\begin{equation*}
\left\vert \varphi (\lambda )-x_{m}\right\vert =\left\vert x_{\left\vert 
\mathbb{N}\cap \lambda \right\vert }-x_{m}\right\vert <\varepsilon
\end{equation*}%
and taking the $\Lambda $-limit, we get the conclusion: 
\begin{equation*}
\varepsilon >\lim_{\lambda \uparrow \Lambda }\ \left\vert \varphi (\lambda
)-x_{m}\right\vert =\left\vert \lim_{\lambda \uparrow \Lambda }\ \varphi
(\lambda )-x_{m}\right\vert =\left\vert L-x_{m}\right\vert .
\end{equation*}

$\square $

\bigskip

From now on, the set $\mathbb{R}$ of standard numbers will be identified
with the set of real numbers, namely the real number will be considered
"special" points on the Euclidean line.

Given a net $\varphi :\Lambda \rightarrow \mathbb{R}$, since $\Lambda $ is a
directed set, also the Cauchy limit is well defined: 
\begin{equation}
L=\lim_{\lambda \rightarrow \Lambda }\varphi (\lambda )\Leftrightarrow
\forall \varepsilon \in \mathbb{R}^{+},\exists \lambda _{0}\in \mathfrak{L}%
,\forall \lambda \supset \lambda _{0},|\varphi (\lambda )-L|\ \leq
\varepsilon  \label{limC}
\end{equation}

Notice that in order to distinguish the Cauchy limit (\ref{limC}) from the $%
\Lambda $-limit, we have used the symbols "$\lambda \rightarrow \Lambda $"
and "$\lambda \uparrow \Lambda $" respectively.

The standard part of a number is related to the Cauchy notion of limit. If a
real net $x_{\lambda }$ admits the Cauchy limit, the relation with the $%
\Lambda $-limit is given by the following identity:%
\begin{equation}
\lim_{\lambda \rightarrow \Lambda }\ x_{\lambda }=st\left( \lim_{\lambda
\uparrow \Lambda }\ x_{\lambda }\right)  \label{bn}
\end{equation}

Another important relation between the two limits is the following:

\begin{proposition}
\label{lim}If 
\begin{equation*}
\lim_{\lambda \uparrow \Lambda }\ x_{\lambda }=\xi \in \mathbb{E}_{\kappa }
\end{equation*}%
and $\xi $ is bounded, then there exist a sequence $\lambda _{n}\in 
\mathfrak{L}$ such that%
\begin{equation*}
\lim_{n\rightarrow \infty }\ x_{\lambda _{n}}=st(\xi ).
\end{equation*}
\end{proposition}

\textbf{Proof: }Set $x_{0}=st(\xi )$ and for every $n\in \mathbb{N}$, take $%
\lambda _{n}$ such that $\ x_{\lambda _{n}}\in \left[ x_{0}-1/n,\ x_{0}+1/n%
\right] .$

$\square $

\begin{remark}
As we have remarked in the intruduction, the Centrum Axiom is necessary to
prove Theorem \ref{sat}. Actually it is not difficult to prove that the
Centrum Axiom is equivalent to the following:%
\begin{equation*}
\mathbb{E_{\kappa }\ }\text{contains a subfield isomorphic to the field of
real numbers.}
\end{equation*}
\end{remark}

\bigskip

The notion of "standard entity" can be extended from numbers (i.e. the real
numbers) to other elements of the universe by the following definition:

\begin{definition}
An element $E\in V(\mathbb{R})$ is called \textbf{standard }and\textbf{\ }$V(%
\mathbb{R})$ is called \textbf{standard universe;} $V(\mathbb{E}_{\kappa })$
is called \textbf{Euclidean universe}.
\end{definition}

Notice that 
\begin{equation*}
V(\mathbb{R})\subset \Lambda \subset V(\mathbb{E}_{\kappa }).
\end{equation*}%
Also the second inclusion is strict since $V(\mathbb{E}_{\kappa })$ contains
sets of inaccessible cardinality such as $\mathbb{E}_{\kappa }$.

\section{The Euclidean universe}

In this section we will inestigate the structure of the Euclidean universe $%
V(\mathbb{E}_{\kappa }).$

\subsection{$\Lambda $-limit of sets\label{ls}}

By Def. \ref{lim+}, the $\Lambda $-limit has been defined for every net $%
\varphi \in \mathfrak{F}\left( \mathfrak{L},\mathbb{E}_{\kappa }\right) ;$
next we will extend this notion to the nets of sets in%
\begin{equation*}
\mathfrak{F}\left( \mathfrak{L},V_{n}(\mathbb{E}_{\kappa })\right) =\left\{
\Phi \in V_{n}\left( \mathbb{E}_{\kappa }\right) ^{\mathfrak{L}}\ |\ \exists
A\in \Lambda ,\forall \lambda \in \mathfrak{L},\ \Phi (\lambda )=\Phi
(\lambda \cap A)\right\}
\end{equation*}%
for every $n\in \mathbb{N}$. In the following, in order to simplify the
notation, a net of sets $\Phi \in \mathfrak{F}\left( \mathfrak{L},V_{n}(%
\mathbb{E}_{\kappa })\right) $ will be denoted by $\left\{ E_{\lambda
}\right\} $ where $E_{\lambda }=\Phi \left( \lambda \right) .$

We define the $\Lambda $-limit of sets by induction over $n.$ If $n=0,$ $%
\lim_{\lambda \uparrow \Lambda }$ $\Phi \left( \lambda \right) $ is a net of
numbers defined by Def. \ref{lim+}; if $n>0,$ we set%
\begin{equation}
E_{\Lambda }=\lim_{\lambda \uparrow \Lambda }\ E_{\lambda }:=\left\{
\lim_{\lambda \uparrow \Lambda }\ \Psi (\lambda )\ |\ \Psi \in \mathfrak{F}%
\left( \mathfrak{L},V_{n-1}(\mathbb{E}_{\kappa })\right) ,\forall \lambda
\in \mathfrak{L}:\Psi (\lambda )\in E_{\lambda }\right\} .  \label{vivaldi}
\end{equation}%
Clearly, by (\ref{l1}), $E_{\Lambda }\in V(\mathbb{E}_{\kappa }).$

\begin{definition}
A set $E$ obtained as $\Lambda $-limit of a net of sets $\left\{ E_{\lambda
}\right\} \in \mathfrak{F}\left( \mathfrak{L},V_{n}(\mathbb{E}_{\kappa
})\right) $ is called \textbf{internal}. If not it, is called \textbf{%
external}.
\end{definition}

For example the set $\mathbb{R}$ is external.

If $C_{A}\left( \lambda \right) =A\in \Lambda _{S}$ is a constant net, we
set 
\begin{equation}
A^{\ast }:=\lim_{\lambda \uparrow \Lambda }\ C_{A}\left( \lambda \right)
=\left\{ \lim_{\lambda \uparrow \Lambda }\ \Psi (\lambda )\ |\ \Psi \in 
\mathfrak{F}\left( \mathfrak{L},V_{n-1}(\mathbb{E}_{\kappa })\right)
,\forall \lambda \in \mathfrak{L}:\Psi (\lambda )\in A\right\} ;
\label{Astar}
\end{equation}%
then, if $A\in V_{n}(\mathbb{E}_{\kappa }),$ also $A^{\ast }\in V_{n}(%
\mathbb{E}_{\kappa }).$ This definition extends (\ref{star-}) to all the
elements of $\Lambda =\Lambda _{S}\cup \mathbb{E}_{\kappa }$. $A^{\ast }$
will be called the $\ast $\textbf{-transform} of $A.$

The $\ast $-transform\textbf{\ }allows to build a family $\left\{ \mathbb{E}%
_{j}\right\} _{j\in \mathbf{Ord}}$ of subsets of $\mathbb{E}_{\kappa }$ as
follows:

\begin{itemize}
\item $\mathbb{E}_{0}=\mathbb{R};$

\item $\mathbb{E}_{j+1}=\mathbb{E}_{j}^{\ast };$

\item if $j\leq \kappa $ is a limit ordinal, then $\mathbb{E}%
_{j}=\dbigcup\limits_{k<j}\mathbb{E}_{k}.$
\end{itemize}

\subsection{$\Lambda $-limit of functions}

Since in set theory a function $f$ can be identified with its graph $\Gamma
_{f},$%
\begin{equation*}
f_{\Lambda }:=\lim_{\lambda \uparrow \Lambda }~f_{\lambda }
\end{equation*}%
is well defined. However, it is not immediate to see that $f_{\Lambda }$ is
function. For this reason, we will analyze this situation explicitly.

\begin{theorem}
Given a net of functions $\left\{ f_{\lambda }\right\} $ 
\begin{equation*}
f_{\lambda }:A_{\lambda }\rightarrow B_{\lambda },\ \ A_{\lambda
},B_{\lambda }\in V_{n}(\mathbb{E}_{\kappa }),
\end{equation*}%
then $f_{\Lambda }:A_{\Lambda }\rightarrow B_{\Lambda }$ defined by 
\begin{equation}
f_{\Lambda }\left( \lim_{\lambda \uparrow \Lambda }\ x_{\lambda }\right)
=\lim_{\lambda \uparrow \Lambda }~f_{\lambda }\left( x_{\lambda }\right) ;
\label{42}
\end{equation}%
is a function and we have that%
\begin{equation*}
\Gamma _{f_{\Lambda }}=\left( \Gamma _{f}\right) _{\Lambda }
\end{equation*}
\end{theorem}

\textbf{Proof}: First, we will prove that (\ref{42}) is a good definition,
namely that $f_{\Lambda }(\xi )$ does not depend on the net $x_{\lambda }$
which defines $\xi .$ We set%
\begin{equation*}
\xi =\lim_{\lambda \uparrow \Lambda }\ x_{\lambda }=\lim_{\lambda \uparrow
\Lambda }\ y_{\lambda }
\end{equation*}%
and we have to prove that%
\begin{equation*}
\lim_{\lambda \uparrow \Lambda }\ f\left( x_{\lambda }\right) =\lim_{\lambda
\uparrow \Lambda }\ f\left( y_{\lambda }\right) .
\end{equation*}%
We take 
\begin{equation*}
\chi \left( \lambda \right) =\left\{ 
\begin{array}{cc}
1 & if\ \ x_{\lambda }=y_{\lambda } \\ 
0 & if\ \ x_{\lambda }\neq y_{\lambda }%
\end{array}%
\right. .
\end{equation*}%
Hence $\forall \lambda $, $\chi \left( \lambda \right) +\left( x_{\lambda
}-y_{\lambda }\right) \neq 0$ and so%
\begin{equation*}
\lim_{\lambda \uparrow \Lambda }\ \left[ \chi \left( \lambda \right) +\left(
x_{\lambda }-y_{\lambda }\right) \right] \neq 0,
\end{equation*}%
then,%
\begin{eqnarray*}
\lim_{\lambda \uparrow \Lambda }\chi \left( \lambda \right) &=&\lim_{\lambda
\uparrow \Lambda }\chi \left( \lambda \right) +\lim_{\lambda \uparrow
\Lambda }\ x_{\lambda }-\lim_{\lambda \uparrow \Lambda }\ y_{\lambda } \\
&=&\lim_{\lambda \uparrow \Lambda }\ \left[ \chi \left( \lambda \right)
+\left( x_{\lambda }-y_{\lambda }\right) \right] \neq 0.
\end{eqnarray*}%
Moreover, we have that 
\begin{equation*}
\forall \lambda ,\ \chi \left( \lambda \right) \cdot \left[ f_{\lambda
}\left( x_{\lambda }\right) -\ f_{\lambda }\left( y_{\lambda }\right) \right]
=0;
\end{equation*}%
then%
\begin{equation*}
0=\lim_{\lambda \uparrow \Lambda }\left( \chi \left( \lambda \right) \cdot 
\left[ f\left( x_{\lambda }\right) -\ f\left( y_{\lambda }\right) \right]
\right) =\lim_{\lambda \uparrow \Lambda }\chi \left( \lambda \right) \cdot
\lim_{\lambda \uparrow \Lambda }\left[ f_{\lambda }\left( x_{\lambda
}\right) -f_{\lambda }\left( y_{\lambda }\right) \right] .
\end{equation*}%
Since $\lim_{\lambda \uparrow \Lambda }\chi \left( \lambda \right) \neq 0,$
we have that%
\begin{equation*}
0=\lim_{\lambda \uparrow \Lambda }\left[ f_{\lambda }\left( x_{\lambda
}\right) -\ f_{\lambda }\left( y_{\lambda }\right) \right] =\lim_{\lambda
\uparrow \Lambda }\ f_{\lambda }\left( x_{\lambda }\right) -\lim_{\lambda
\uparrow \Lambda }\ f_{\lambda }\left( y_{\lambda }\right) .
\end{equation*}%
Finally, it is immediate to check that $f^{\ast }$ is the graph of the
function (\ref{42}), in fact%
\begin{eqnarray*}
\left( \Gamma _{f}\right) _{\Lambda } &=&\left\{ \lim_{\lambda \uparrow
\Lambda }\left( x_{\lambda },f_{\lambda }\left( x_{\lambda }\right) \right)
\ |\ \forall \lambda ,\ \left( x_{\lambda },f_{\lambda }\left( x_{\lambda
}\right) \right) \in \Gamma _{f_{\lambda }}\right\} \\
&=&\left\{ \left( \lim_{\lambda \uparrow \Lambda }\ x_{\lambda
},\lim_{\lambda \uparrow \Lambda }f_{\lambda }\left( x_{\lambda }\right)
\right) \ |\ \forall \lambda ,\ x_{\lambda }=f_{\lambda }\left( x_{\lambda
}\right) \right\} \\
&=&\left\{ \left( \xi ,\ f_{\Lambda }\left( \xi \right) \right) \ |\ \xi
=f_{\Lambda }\left( \xi \right) \right\} =\Gamma _{f_{\Lambda }}.
\end{eqnarray*}

$\square $

\bigskip

\begin{definition}
A function $f$ obtained as $\Lambda $-limit of a net of functions $%
f_{\lambda }$ is called \textbf{internal}. Otherwise is called \textbf{%
external}.
\end{definition}

If $\left\{ f\right\} $ is a constant net, we set 
\begin{equation}
f^{\ast }=\lim_{\lambda \uparrow \Lambda }\ f;
\end{equation}%
then, if $f:A\rightarrow B,$ then $f^{\ast }:A^{\ast }\rightarrow B^{\ast }$
and $f^{\ast }$ will be called the $\ast $\textbf{-transform} of $f.$

\subsection{Hyperfinite sets}

Another fundamental notion in Euclidean calculus is the following:

\begin{definition}
We say that a set $F\in V\left( \mathbb{E}_{\kappa }\right) $ is \textbf{%
hyperfinite} if there is a net $\left\{ F_{\lambda }\right\} _{\lambda \in
\Lambda }$ of finite sets such that 
\begin{equation*}
F=\lim_{\lambda \uparrow \Lambda }~F_{\lambda }=\left\{ \lim_{\lambda
\uparrow \Lambda }\ x_{\lambda }\ |\ x_{\lambda }\in F_{\lambda }\right\}
\end{equation*}
\end{definition}

The hyperfinite sets share many properties of finite sets. For example, a
hyperfinite set $F\subset \mathbb{E}_{\kappa }$ has a maximum $x_{M}$ and a
minimum $x_{m}$ respectively given by%
\begin{equation*}
x_{M}=\lim_{\lambda \uparrow \Lambda }\max F_{\lambda };\ \
x_{m}=\lim_{\lambda \uparrow \Lambda }\min F_{\lambda }
\end{equation*}

Moreover, it is possible to "add" the elements of an hyperfinite set of
numbers. If $F$ is an hyperfinite set of numbers, the \textbf{hyperfinite sum%
} of the elements of $F$ is defined as follows: 
\begin{equation*}
\sum_{x\in F}x=\ \lim_{\lambda \uparrow \Lambda }\sum_{x\in F_{\lambda }}x.
\end{equation*}

One peculiarity of Euclidean analysis the possibility to associate a unique
hyperfinite set $E^{\circledcirc }$ to any set $E\in V(\mathbb{E}_{\kappa })$
according to the following definition:

\begin{definition}
\label{tondo}Given a set $E\in \Lambda _{S},$ the set 
\begin{equation*}
E^{\circledcirc }:=\lim_{\lambda \uparrow \Lambda }\ \left( E\cap \lambda
\right)
\end{equation*}%
is called \textbf{hyperfinite extension} of $E.$
\end{definition}

If $F=\lim_{\lambda \uparrow \Lambda }~F_{\lambda }$ is a hyperfinite set,
its \textbf{hypercardinality} is given by%
\begin{equation*}
\left\vert F\right\vert ^{\ast }=\lim_{\lambda \uparrow \Lambda }~\left\vert
F_{\lambda }\right\vert
\end{equation*}%
were $\left\vert \cdot \right\vert ^{\ast }$ is the $\ast $-tranform of the
fuction "cardinality" defined on finite sets.

Notice that, by vitue of (\ref{bella}), the hypercardinality of $%
E^{\circledcirc },$ given by%
\begin{equation*}
\left\vert E^{\circledcirc }\right\vert ^{\ast }=\lim_{\lambda \uparrow
\Lambda }\ \left\vert E\cap \lambda \right\vert ,
\end{equation*}%
is the numerosity of $E$ as it has been defined by Axiom \ref{num}.

If we put%
\begin{equation*}
E^{\sigma }=\left\{ x^{\ast }\ |\ x\in E\right\} ,
\end{equation*}%
we can associate the sets $E^{\sigma },$ $E^{\circledcirc }$ and $E^{\ast }$
to any set $E\in \Lambda .$ They are ordered as follows:%
\begin{equation*}
E^{\sigma }\subseteq E^{\circledcirc }\subseteq E^{\ast };
\end{equation*}%
in particular, if $E\subseteq \mathbb{R}$, $E^{\sigma }=E.$ The hyperfinite
analysis is very relevant in the applications and the operator "$%
\circledcirc $" plays a special role. You can see some examples of this fact
in section \ref{i}.

\section{Nonstandard Analysis}

Even if some notions and definitions of Nonstandard Analysis have already
been introduced in the previous sections, now we will treat this topic in
details. In this section, we assume that the reader is familiar with the
basic notions on NSA.

\subsection{Nonstandard theories}

In this subsection, we recall the basic notion of Nonstandard Analysis how
have been developed in the "superstructure" approach. Following Keisler (see 
\cite{keisler}), we give the following definition:

\begin{definition}
\label{nu}A nonstandard theory is a triple $\left( V\left( \mathbb{R}\right)
,V\left( \mathbb{R}^{\mathbb{\bullet }}\right) ,\bullet \right) $ such that%
\footnote{%
To be precise, Keisler calls $\left( V\left( \mathbb{R}\right) ,V\left( 
\mathbb{R}^{\mathbb{\bullet }}\right) ,\bullet \right) $ \textit{nonstandard
universe} while we use the espression \textit{nonstandard universe} to
denote the set $V\left( \mathbb{R}^{\mathbb{\bullet }}\right) $.}

\begin{itemize}
\item $V\left( \mathbb{R}\right) $ is a superstructure over $\mathbb{R}$
called standard universe;

\item $\mathbb{R}^{\mathbb{\bullet }}$ is a set such that $\mathbb{R\subset R%
}^{\mathbb{\bullet }}$ which is called field of the ($\bullet $)-hyperreal
numbers;

\item $V\left( \mathbb{R}^{\mathbb{\bullet }}\right) $ is a superstructure
over $\mathbb{R}^{\mathbb{\bullet }}$ called nonstandard universe;

\item the map 
\begin{equation*}
\bullet ~:V\left( \mathbb{R}\right) \rightarrow V\left( \mathbb{R}^{\mathbb{%
\bullet }}\right)
\end{equation*}%
satisfies the Leibniz principle and 
\begin{equation}
\forall r\in \mathbb{R},r=r^{\mathbb{\bullet }}.  \label{ul}
\end{equation}
\end{itemize}
\end{definition}

We recall the notion of Leibniz (or transfer) Principle. It is well known
that the map $\bullet $ transforms any elementary sentence $%
P(a_{1},a_{2},...,a_{n})$ to a elementary sentence $P(a_{1}^{\mathbb{\bullet 
}},a_{2}^{\mathbb{\bullet }},...,a_{n}^{\mathbb{\bullet }})$ in $V\left( 
\mathbb{R}^{\mathbb{\bullet }}\right) $ where $a_{1},a_{2},...,a_{n}$ are
constants in $V\left( \mathbb{R}\right) .$ The adjective elementary refers
to the fact that the quantifiers in elementary sentences are of the form ($%
\forall x\in y$) or ($\exists x\in y$) where $x$ is a variable and $y$ is a
constant or a variable. The Leibniz principle states that $%
P(a_{1},a_{2},...,a_{n})$ is true if an only if $P(a_{1}^{\mathbb{\bullet }%
},a_{2}^{\mathbb{\bullet }},...,a_{n}^{\mathbb{\bullet }})$ is true. For
details, see e.g. \cite{BDN2018} or \cite{keisler}.

\begin{definition}
Given two sets $\mathbb{A}$ and $\mathbb{S}$, a superstructure embedding is
a triple $\left( V\left( \mathbb{A}\right) ,V\left( \mathbb{S}\right)
,\bullet \right) $ where$\ \bullet :V\left( \mathbb{A}\right) \rightarrow
V\left( \mathbb{S}\right) $ is a injective map such that%
\begin{equation*}
\mathbb{A^{\mathbb{\bullet }}=S}
\end{equation*}%
\begin{equation*}
\forall x,y\in V\left( \mathbb{A}\right) ,\ x\in y\Leftrightarrow x^{\mathbb{%
\bullet }}\in y^{\mathbb{\bullet }}
\end{equation*}
\end{definition}

The following fact is well known:

\begin{theorem}
\label{logic}If $\left( V\left( \mathbb{A}\right) ,V\left( \mathbb{S}\right)
,\bullet \right) $ is superstructure embedding then the map $\bullet $
satisfies the Leibniz principle.
\end{theorem}

\textbf{Proof}: This result can be proved by induction over the complexity
of the sentences; see e.g. \cite{BDN2018} Th. 5.8. or \cite{keisler}.

$\square $\bigskip

By the above theorem, we get the following result:

\begin{corollary}
\label{logic2}If $\left( V\left( \mathbb{R}\right) ,V\left( \mathbb{S}%
\right) ,\bullet \right) $ is a superstructure embedding such that $\mathbb{%
R\neq S}$,\ then, $\left( V\left( \mathbb{R}\right) ,V\left( \mathbb{R^{%
\mathbb{\bullet }}}\right) ,\bullet \right) $ is a nonstandard theory with $%
\mathbb{R^{\mathbb{\bullet }}=S}$.
\end{corollary}

An isomorphism beetwen nonstandard theories is defined as follows:

\begin{definition}
Two nonstandard theories $\left( V\left( \mathbb{R}\right) ,V\left( \mathbb{R%
}^{\mathbb{\bullet }}\right) ,\bullet \right) $ and $\left( V\left( \mathbb{R%
}\right) ,V\left( \mathbb{R}^{\star }\right) ,\star \right) $ are isomorphic
if there is a map%
\begin{equation*}
h:V\left( \mathbb{R}^{\mathbb{\bullet }}\right) \rightarrow V\left( \mathbb{R%
}^{\star }\right)
\end{equation*}%
such that

\begin{enumerate}
\item $\forall r\in \mathbb{R},$ $h(r)=r;$

\item $h$ maps $\mathbb{R}^{\mathbb{\bullet }}$ one to one onto $\mathbb{R}%
^{\star };$

\item for each $A\in V\left( \mathbb{R}^{\mathbb{\bullet }}\right)
\backslash \mathbb{R}^{\mathbb{\bullet }},$%
\begin{equation*}
h(A)=\left\{ h(a)\ |\ a\in A\right\}
\end{equation*}

\item for each $A\in V\left( \mathbb{R}\right) ,$ $h(A^{\mathbb{\bullet }%
})=A^{\star }.$
\end{enumerate}
\end{definition}

\begin{definition}
A nonstandard theory $\left( V\left( \mathbb{R}\right) ,V\left( \mathbb{R}^{%
\mathbb{\bullet }}\right) ,\bullet \right) $ is called \textbf{saturated} if
any family of sets $\mathfrak{S}\in V_{n}\left( \mathbb{R}\right) ^{\mathbb{%
\bullet }}$ with cardinality smaller than $\mathbb{R}^{\mathbb{\bullet }}$
and with the finite intersection property has non empty intersection; namely
if 
\begin{equation*}
S_{1}\cap ...\cap S_{n}\neq \varnothing ,\ \ S_{i}\in \mathfrak{S};\
\left\vert \mathfrak{S}\right\vert <\left\vert \mathbb{R}^{\mathbb{\bullet }%
}\right\vert
\end{equation*}%
then%
\begin{equation*}
\dbigcap \mathfrak{S}\neq \varnothing .
\end{equation*}
\end{definition}

Among all the nonstandard theories there is a privileged one which is unique
up to isomorphisms.

\begin{theorem}
\label{ki}A saturated nonstandard theory $\left( V\left( \mathbb{R}\right)
,V\left( \mathbb{R}^{\bullet }\right) ,\bullet \right) $ with 
\begin{equation*}
\left\vert \mathbb{R}^{\mathbb{\bullet }}\right\vert =\kappa
\end{equation*}%
is unique up to isomorphism. In this case, $V\left( \mathbb{R}^{\bullet
}\right) $ will be called Keisler universe.
\end{theorem}

\textbf{Proof}: See \cite{keisler}.

$\square $

\subsection{The Normal Universe}

According to the theory of the previous section we give following

\begin{definition}
If "$\ast $" is the map defined by (\ref{Astar}), $V\left( \mathbb{R}^{\ast
}\right) $ will be called \textbf{normal universe}; $\mathbb{R}^{\ast }$%
\textbf{\ }will be called \textbf{normal Euclidean field }and from now on,
it will be simply denoted by $\mathbb{E}$.
\end{definition}

\begin{theorem}
\label{satur}The triple $\left( V\left( \mathbb{R}\right) ,V\left( \mathbb{E}%
\right) ,\ast \right) $ is a nonstandard theory and $V\left( \mathbb{E}%
\right) $ is a Keisler universe.
\end{theorem}

\textbf{Proof}: By the definition of $V\left( \mathbb{E}\right) ,$ it
follows that $\left( V\left( \mathbb{R}\right) ,V\left( \mathbb{E}\right)
,\ast \right) $ is a superstructure embedding. Then by Cor. \ref{logic2}, we
have to prove that $\left\vert \mathbb{E}\right\vert =\kappa $ and that $%
\left( V\left( \mathbb{R}\right) ,V\left( \mathbb{E}\right) ,\ast \right) $
is saturated. Since 
\begin{equation*}
\mathbb{E=R}^{\ast }=\left\{ \lim_{\lambda \rightarrow \Lambda }\ x_{\lambda
}\ |\ x_{\lambda }\in \mathfrak{F}\left( \mathfrak{L},\mathbb{R}\right)
\right\}
\end{equation*}%
we have that 
\begin{equation*}
\left\vert \mathbb{E}\right\vert \leq \left\vert \mathfrak{F}\left( 
\mathfrak{L},\mathbb{R}\right) \right\vert
\end{equation*}%
Moreover%
\begin{equation*}
\mathfrak{F}\left( \mathfrak{L},\mathbb{R}\right) =\left\{ \varphi \in 
\mathbb{R}^{\mathfrak{L}}\ |\ \exists A\in \Lambda ,\forall \lambda \in 
\mathfrak{L},\ \varphi (\lambda )=\varphi (\lambda \cap A)\right\}
\end{equation*}%
and hence%
\begin{equation*}
\left\vert \mathfrak{F}\left( \mathfrak{L},\mathbb{R}\right) \right\vert
=\left\vert \dbigcup\limits_{A\in \Lambda }\mathbb{R}^{A}\right\vert
\end{equation*}%
and since $\left\vert \mathbb{R}^{A}\right\vert <\kappa $ and $\left\vert
\Lambda \right\vert =\kappa $, we have that%
\begin{equation*}
\left\vert \mathbb{E}\right\vert \leq \left\vert \dbigcup\limits_{A\in
\Lambda }\mathbb{R}^{A}\right\vert =\kappa
\end{equation*}%
Also we have that 
\begin{equation*}
\left\vert \mathbb{E}\right\vert \geq \left\vert \mathfrak{N}\right\vert
\geq \left\vert \mathbf{Ord}\right\vert =\kappa .
\end{equation*}%
Then $\left\vert \mathbb{E}\right\vert =\kappa $. It remains to prove that
it is saturated.

If $\mathfrak{S}\in V_{n}\left( \mathbb{R}\right) ^{\mathbb{\ast }},$ then 
\begin{equation*}
\mathfrak{S}=\left\{ E_{\mu }\ |\ \mu \in H\right\}
\end{equation*}%
where $H$ is a set of indices with $\left\vert H\right\vert <\kappa $. Since 
$\left\vert H\right\vert <\kappa ,$ it is not restrictive to assume that $%
H\subset \mathfrak{L}$. For every $\mu \in H$, let $\varphi _{\mu }(\lambda
) $ be a net such that 
\begin{equation*}
\lim_{\lambda \uparrow \Lambda }\ \varphi _{\mu }(\lambda )=E_{\mu }
\end{equation*}%
For any fixed $\lambda $, pick an element%
\begin{equation*}
\psi (\lambda )\in \dbigcap_{\mu \subseteq \lambda }\varphi _{\mu }(\lambda )
\end{equation*}%
if this intersection is nonempty. Otherwise, pick 
\begin{equation*}
\psi (\lambda )\in \dbigcap_{\mu \subseteq \lambda \backslash \left\{
x_{1}\right\} }\varphi _{\mu }(\lambda );\ x_{1}\in \lambda
\end{equation*}%
if this intersection is nonempty, and continue in this manner until the
element $\psi (\lambda )$ is defined. In case that this intersection is
always empty, we set $\psi (\lambda )=\varnothing .$ As a consequence of
this definition, the following property holds:%
\begin{equation}
\dbigcap_{\mu \subseteq \tau }\varphi _{\mu }(\lambda )\neq \varnothing
\Rightarrow \forall \lambda \supseteq \tau ,\ \psi (\lambda )\in
\dbigcap_{\mu \subseteq \tau }\varphi _{\mu }(\lambda )  \label{it}
\end{equation}%
Now let $\tau \in H$ be fixed. By the finite intersection property,%
\begin{equation*}
\varnothing \neq \dbigcap_{\mu \subseteq \tau }E_{\mu }=\dbigcap_{\mu
\subseteq \tau }\lim_{\lambda \uparrow \Lambda }\ \varphi _{\mu }(\lambda
)=\lim_{\lambda \uparrow \Lambda }\left( \dbigcap_{\mu \subseteq \tau
}\varphi _{\mu }(\lambda )\right)
\end{equation*}%
Then, there exists a set $Q\in \mathcal{U}$ ($\mathcal{U}$ is defined in the
proof of Th. \ref{LL}) such that, 
\begin{equation*}
\forall \lambda \in Q,\ \dbigcap_{\mu \subseteq \tau }\varphi _{\mu
}(\lambda )\neq \varnothing ;
\end{equation*}%
and hence, by (\ref{it}), it follows that 
\begin{equation*}
\forall \lambda \in Q,\ \psi (\lambda )\in \dbigcap_{\mu \subseteq \tau
}\varphi _{\mu }(\lambda )\neq \varnothing
\end{equation*}%
and taking the $\Lambda $-limit 
\begin{equation*}
\lim_{\lambda \uparrow \Lambda }\ \psi (\lambda )\in \dbigcap_{\mu \subseteq
\tau }\lim_{\lambda \uparrow \Lambda }\ \varphi _{\mu }(\lambda
)=\dbigcap_{\mu \subseteq \tau }E_{\mu }
\end{equation*}%
and in particular, 
\begin{equation*}
\lim_{\lambda \uparrow \Lambda }\ \psi (\lambda )\in E_{\tau }
\end{equation*}

As this holds for every $\tau \in H$, we conclude that%
\begin{equation*}
\lim_{\lambda \uparrow \Lambda }\ \psi (\lambda )\in \dbigcap_{\tau \in
H}E_{\tau }=\dbigcap \mathfrak{S.}
\end{equation*}%
$\square $\bigskip

$V\left( \mathbb{E}\right) $ is not the only Keisler universe contained in $%
V\left( \mathbb{E}_{\kappa }\right) $. For example, using the notations at
the end of section \ref{ls}, we have that also 
\begin{equation*}
\left( V\left( \mathbb{R}\right) ,V\left( \mathbb{E}_{2}\right) ,\ast \ast
\right) ,\ \left( V\left( \mathbb{R}\right) ,V\left( \mathbb{E}_{3}\right)
,\ast \ast \ast \right) ,\ \ etc.
\end{equation*}%
are saturated nonstandard theories and hence we have infinite Keisler
universes included in $V\left( \mathbb{E}_{\kappa }\right) $. More in
general, for $0<j\leq \kappa $, it is possible to define a Keisler universe%
\begin{equation*}
\left( V\left( \mathbb{R}\right) ,V\left( \mathbb{E}_{j}\right) ,\ast
^{j}\right)
\end{equation*}%
where the map $\ast ^{j}$ is defined $\forall A\in V\left( \mathbb{R}\right) 
$ as follows:

\begin{itemize}
\item if $j=1,$ $A^{\ast ^{1}}=A^{\ast }\subset \mathbb{E}_{1}(=\mathbb{E})$

\item if $j=k+1,$ $A^{\ast ^{j+1}}=\left( A^{\ast ^{j}}\right) ^{\ast
}\subset \mathbb{E}_{j+1}$

\item if $k\ $is a limit ordinal, and $A\in V\left( \mathbb{E}_{k}\right) ,$
then $A\in V\left( \mathbb{E}_{j}\right) $ for some $j<k$ and hence $A^{\ast
^{k}}=A^{\ast ^{j}}.$
\end{itemize}

In particular we have that all the fields $\mathbb{E}_{j},$ $j\leq \kappa $
are isomorphic; for $j<\kappa ,$ the map%
\begin{equation*}
\ast ~:\mathbb{E}_{j}\rightarrow \mathbb{E}_{j+1}
\end{equation*}%
is a field homomorphism and if $k$ is a limit ordinal, the map%
\begin{equation*}
\ast ~:\mathbb{E}_{k}\rightarrow \mathbb{E}_{k}
\end{equation*}%
is a field isomorphism. In any case, the only fixed points of $\ast $ are
the real numbers. The spaces $\mathbb{E}_{j}$'s differ from each other by
the way they are embedded in $\mathbb{A}$.

Moreover in a Euclidean universe there are other interesting superstructure
embeddings which can be useful in some application. For example, 
\begin{equation*}
\left( V\left( \mathbb{E}_{j}\right) ,V\left( \mathbb{E}_{j+1}\right) ,\ast
\right)
\end{equation*}%
is a superstructure embedding; however is not a nonstandard theory, since $%
\mathbb{E}_{j}\neq \mathbb{R}$; moreover $\left( V\left( \mathbb{E}%
_{j}\right) ,V\left( \mathbb{E}_{j+1}\right) ,\ast \right) $ violates an
other important request of Def. \ref{nu}, namely 
\begin{equation*}
\exists \xi \in \mathbb{E}_{j},\ \xi \neq \xi ^{\ast }
\end{equation*}%
in contrast with (\ref{ul}). Nevertheless, by Th. \ref{logic}, $\left(
V\left( \mathbb{E}_{j}\right) ,V\left( \mathbb{E}_{j+1}\right) ,\ast \right) 
$ satisfies the Leibniz principle. Then the following fact follows
straightforwardly:

\begin{theorem}
$\left( V\left( \mathbb{E}_{\kappa }\right) ,V\left( \mathbb{E}_{\kappa
}\right) ,\ast \right) $ is a superstructure embedding that satisfies the
Leibniz principle.
\end{theorem}

\subsection{The $\protect\alpha $-theory\label{AT}}

The $\alpha $-theory has been introduced in \cite{BDN03} and it represents
an elementary approach to nonstandard analysis particularly suitable for
some application; see e.g \cite{BDN2018}, \cite{BGG}, \cite{BC2021} and
references. Actually, the Euclidean universe contains many non-standard
universes which can be easily defined and, therefore, are more suitable for
the elementary applications of practitioners. The $\alpha $-theory is one of
them and it can be constructed using the notion of $\alpha $-limit where $%
\alpha =\mathfrak{num}\left( \mathbb{N}^{+}\right) $, $\mathbb{N}^{+}=%
\mathbb{N}\backslash \left\{ 0\right\} $.

\begin{definition}
Given a sequence $\varphi :$ $\mathbb{N}\rightarrow V_{n}(\mathbb{E}_{\kappa
}),$ we set%
\begin{equation*}
\lim_{n\uparrow \alpha }\ \varphi (n):=\lim_{\lambda \uparrow \Lambda }\
\varphi \left( \left\vert \lambda \cap \mathbb{N}^{+}\right\vert \right) .
\end{equation*}
\end{definition}

If $C_{b}(n)$ is the constant sequence with value $b\in V(\mathbb{E}_{\kappa
}),$ we set%
\begin{equation*}
b^{\ast _{\alpha }}=\lim_{n\uparrow \alpha }\ C_{b}(n)
\end{equation*}%
Then, by Th. \ref{logic}, it follows that%
\begin{equation*}
\left( V\left( \mathbb{R}\right) ,V\left( \mathbb{R}^{\ast _{\alpha
}}\right) ,\ast _{\alpha }\right)
\end{equation*}%
is a nonstandard theory.

\begin{definition}
The nonstandard theory $\left( V\left( \mathbb{R}\right) ,V\left( \mathbb{R}%
^{\ast _{\alpha }}\right) ,\ast _{\alpha }\right) $ is called $\alpha $%
-theory.
\end{definition}

If 
\begin{equation*}
i:\mathbb{N}\rightarrow \mathbb{E}_{\kappa };\ i(n)=n
\end{equation*}%
then taking the $\alpha $-limit we get that 
\begin{equation*}
\lim_{n\uparrow \alpha }\ i(n)=\lim_{\lambda \uparrow \Lambda }\ i\left(
\left\vert \lambda \cap \mathbb{N}^{+}\right\vert \right) =\mathfrak{num}(%
\mathbb{N}^{+})=\alpha
\end{equation*}%
hence the set of ($\ast _{\alpha }$)-hyperreal numbers $\mathbb{R}^{\ast
_{\alpha }}$ can be characterized as follows:%
\begin{equation*}
\mathbb{R}^{\ast _{\alpha }}=\left\{ \lim_{n\uparrow \alpha }\ \varphi (n)\
|\ \varphi \,:\mathbb{N}\rightarrow \mathbb{R}\right\}
\end{equation*}%
namely, every ($\ast _{\alpha }$)-hyperreal numbers is the $\alpha $-limit
of real sequence. For example, we have that $\omega $ is a ($\ast _{\alpha }$%
)-hyperreal number since 
\begin{eqnarray*}
\omega &=&\mathfrak{num}(\mathbb{N})=\mathfrak{num}(\mathbb{N}^{+})+%
\mathfrak{num}(\left\{ 0\right\} ) \\
&=&\alpha +1=\lim_{n\uparrow \alpha }\ (n+1).
\end{eqnarray*}

As we will see in section \ref{SE}, the construction of a model of the
theory, there is an ultrafilter $\mathcal{U}$ which plays a central role
(see definition (\ref{uf})). By choosing $\mathcal{U}$ in a suitable way,
then we get the following result:

\begin{theorem}
\label{anum}It is compatible with axioms 1-3 that the number $\alpha $
satisfies the following properties:

\begin{itemize}
\item \textsc{Divisibility Property} : For every $k\in \mathbb{N}$, the
number $\alpha $ is a multiple of $k$ and the numerosity of the set of
multiples of k:%
\begin{equation*}
\mathfrak{num}(\{k,2k,3k,...,nk,...\})=\frac{\alpha }{k}
\end{equation*}

\item \textsc{Root Property}: For every $k\in \mathbb{N}$, the number $%
\alpha $ is a k-th power and the numerosity of the set of k-th powers:%
\begin{equation*}
\mathfrak{num}(\{1^{k},2^{k},3^{k},...,n^{k},...\})=\sqrt[k]{\alpha }
\end{equation*}

\item \textsc{Power Property}: If we set $\wp _{fin}(A)=\{F\in \wp (A)\ |\
F\ $is a finite set$\},$ then 
\begin{equation*}
\mathfrak{num}(\wp _{fin}(\mathbb{N}^{+}))=2^{\alpha }
\end{equation*}

\item \textsc{Integer numbers Property}: 
\begin{equation*}
\mathfrak{num}(\mathbb{Z})=2\alpha +1
\end{equation*}

\item \textsc{Rational numbers Property}: For every $q\in \mathbb{Q}$,%
\begin{equation*}
\mathfrak{num}(\left( q,q+1\right] \cap \mathbb{Q})=\mathfrak{num}(\left( 0,1%
\right] \cap \mathbb{Q})=\alpha
\end{equation*}%
and 
\begin{equation*}
\mathfrak{num}(\mathbb{Q})=2\alpha ^{2}+1.
\end{equation*}
\end{itemize}
\end{theorem}

\textbf{Proof}: See \cite{BDN2018}, Sections 16.6 and 16.7.

$\square $

\subsection{About the idea of continuum}

The idea of (linear) continuum is described or modeled by the geometric
line. In classical Euclidean geometry, lines and segments are not considered
sets of points; on the contrary, in the last two centuries the reductionist
attitude of modern mathematics has described Euclidean geometry through a
set interpretation. In the last century, the geometric continuum has been
identified with the Dedekind continuum and the geometric line has been
identified with the set of real numbers (once the origin $O$ and a unitary
segment $OU$ have been fixed). Even is this identification, today, is almost
universally accepted, we have seen that also the Euclidean line, as defined
by axiom \ref{1}, has some right to represent the geometric continuum. In
this section we will compare $\mathbb{R}$ and $\mathbb{E\ }$($\mathbb{\cong E%
}_{\kappa }$) with respect to the idea of geometric continuum.

In our naive intuition, we think of a linear continuum as a linearly ordered
set without interruptions, that is, without holes between one part and the
other. Let's make this definition rigorous. Contrary to our intuition, a set
X satisfying the following property%
\begin{equation*}
\forall a,b\in X,\ a<b,\ \exists c\in X,\ a\leq c\leq b
\end{equation*}%
it cannot be considered a continuum: this notion, satisfied for example by
the set of rational numbers, is not a good candidate for a continuum as the
set of rationals is full of holes represented by irrational numbers.

So we are led to discuss the notion of $\kappa $-saturation and to the
Eudoxus Principle:

\begin{definition}
\label{ks}A linearly ordered set $X$ is called $\kappa $-\textbf{saturated}
if it satisfies following property: given two sets $A,B\subset X,$ such that%
\begin{equation}
\left\vert A\right\vert ,\left\vert B\right\vert <\kappa  \label{ab}
\end{equation}%
\begin{equation*}
\forall a\in A,\ \forall b\in B,\ a<b
\end{equation*}%
then $\exists c\in X,$%
\begin{equation*}
\forall a\in A,\ \forall b\in B,\ a\leq c\leq b.
\end{equation*}
\end{definition}

\begin{definition}
(\textbf{Eudoxus Principle}) A linearly ordered Abelian group $\mathbb{F}$
satisfies the Eudoxus Principle if given two sets $A,B\subset \mathbb{F}$
such that%
\begin{equation*}
\forall a\in A,\ \forall b\in B,\ a<b
\end{equation*}%
\begin{equation*}
\forall \varepsilon \in \mathbb{F}^{+},\ \exists a\in A,\ \exists b\in B,\
b-a<\varepsilon .
\end{equation*}%
then $\exists c\in X,$%
\begin{equation*}
\forall a\in A,\ \forall b\in B,\ a\leq c\leq b.
\end{equation*}
\end{definition}

Using these notions, we can characterize $\mathbb{R}$ and $\mathbb{E}$ as
follows:

\begin{theorem}
The field of the real numbers $\mathbb{R}$ is the only field $\mathbb{F}$
such that:

\begin{itemize}
\item (i) satisfies the Eudoxus Principle,

\item (ii) satisfies the Archimedes' Axiom, namely :%
\begin{equation*}
\forall a,b\in \mathbb{F}^{+},\ \exists n\in \mathbb{N},\ na>b
\end{equation*}
\end{itemize}
\end{theorem}

\textbf{Proof}: Well known.

$\square $

\bigskip

The request (ii) is necessary; in fact, for example, the field of rational
functions with a suitable order structure\footnote{%
For example the field of rational fuction $\mathbb{F}$ can be equipped with
an order structure by setting 
\begin{equation*}
\mathbb{F}^{+}=\left\{ \frac{r_{n}x^{n}+r_{n-1}x^{n-1}+....+r_{0}}{%
w_{m}x^{m}+w_{m-1}x^{m-1}+....+w_{0}}\ |\ \frac{r_{n}}{w_{m}}>0\right\} .
\end{equation*}%
} satisfy (i) but not (ii).

\begin{theorem}
The field of the Euclidean numbers $\mathbb{E}$ is the smallest field that:

\begin{itemize}
\item (i) is $\kappa $-saturated,

\item (ii) is a real closed field, namely every polynomial of odd degree has
at least one root.
\end{itemize}
\end{theorem}

\textbf{Proof}: By Th. \ref{satur}, it is easy to check that $\mathbb{E}$ is 
$\kappa $-saturated according to def. \ref{ks}. Moreover, since $\mathbb{E}$
is hyperreal, it is real closed. All the real closed fields of cardinality $%
\kappa $ are isomorphic and $\left\vert \mathbb{E}\right\vert =\kappa $;
hence $\mathbb{E}$ is the smallest of such fields .

$\square $

\bigskip

Notice that the request (ii) is necessary; in fact $\mathbb{Q}^{\ast }$ is a 
$\kappa $-saturated field, but it does not satisfy (i) since the equation $%
x^{3}=2$ does not have any solution in $\mathbb{Q}^{\ast }$. We observe that
the request (ii) fits well the idea of continuity, in fact, a polynomial of
odd degree must take positive and negative values and hence, by continuity,
it must have some $0$'s.

The above discussion suggests the following definitions of continuum:

\begin{definition}
A linearly ordered Abelian group $\mathbb{F}$ is a \textbf{Dedekind continuum%
} if it satisfies the following property: given two sets $A,B\subset X$ such
that%
\begin{equation}
A,B\neq \varnothing  \label{ab'}
\end{equation}%
\begin{equation*}
\forall a\in A,\ \forall b\in B,\ a<b
\end{equation*}%
then $\exists c\in X,$%
\begin{equation*}
\forall a\in A,\ \forall b\in B,\ a\leq c\leq b;
\end{equation*}%
A linearly ordered ordered field $\mathbb{F}$ is an \textbf{absolute
continuum}\footnote{%
This notion of \textit{absolute contiuum} has been introduced by Ehrlich in 
\cite{ehrlich}. However in his definition $\mathbb{F}$ is a class in the
sense of Von Neumann--Bernays--G\"{o}del set theory.} if it is saturated and
real closed.
\end{definition}

With these notions of continuity $\mathbb{R}$ and $\mathbb{E}$ have the
following characterization: $\mathbb{R}$ is the only Dedekind continuum
field; $\mathbb{E}$ is the smallest absolutely continuum field.

\bigskip

\section{Euclidean Calculus}

\bigskip

In this section we pretend to not know the classical calculus and we will
define the basic notion of calculus, derivative and integral, in the most
natural way provided that you are equipped with infinitesimal and infinite
numbers. Hence, these definitions are very similar to those of the XVIII
century. With these definitions, we will discover that every function is
both integrable and it has a left and right derivative. Of course, if a
function is differentiable, the \textit{Euclidean derivative} corresponds to
the usual derivative and if it is Lebesgue-integrable, the \textit{Euclidean
integral} corresponds to the usual Lebesgue integral. We limit this game to
the \textit{normal} functions as defined below.

Of course this game could be extended to other notions and to a larger class
of functions and it might have some interest for the foundations and the
philosophy of Mathematics.

The idea to work directly in a nonstandard unverse is not new; we recall 
\cite{nelson77}, \cite{benci95} and \cite{HR}. This section can be
considered an other experiment in that direction.

\bigskip

\subsection{Normal functions and sets\label{NF}}

\bigskip

In most application, the space $\mathbb{E}_{\kappa }$ and the Euclidean
universe $V\left( \mathbb{E}_{\kappa }\right) $ are too large and hence
might imply useless technicalities. It is more convenient to work in the
normal Euclidean field $\mathbb{E}=\mathbb{R}^{\ast }$ and in the normal
universe $V\left( \mathbb{E}\right) $. So we are lead to the following
definition:

\begin{definition}
A function $f:\mathbb{E\rightarrow E}$ is called \textbf{normal }if%
\begin{equation*}
f=h^{\ast }
\end{equation*}%
where $h\ $is a standard real function, i.e. $f\in \mathbb{R}^{E},\
E\subseteq \mathbb{R}.$
\end{definition}

If $f$ is normal then $\forall x\in \mathbb{R}$, $f\left( x\right) \in 
\mathbb{R}$.

\begin{definition}
A subset $N\subset \mathbb{E}$ is called \textbf{normal }if $N=A^{\ast }$
for some set $A\subset \mathbb{R}$.
\end{definition}

\begin{remark}
In the nonstandard analysis community there is the habit to call standard
both functions and sets of the form $f^{\ast },A^{\ast }$ and functions and
sets in $V(\mathbb{R}).$ Here we call standard the elements of $V(\mathbb{R}%
) $ and normal their counterpart defined as above.
\end{remark}

The usual functions used in the applications of mathematics can be regarded
as normal functions and not as standard functions. The advantadge of this
point of view is that the main notions of infinitesimal analysis can be
defined using the "actual infinitesimal" in a natural way and hence they
assume a different meaning.

\subsection{The notion of derivative}

Since the normal functions are in a biunivocal correspondence with the real
functions, sometimes we will denote both with the same symbol. The same we
will do with the intervals.

Now let us introduce some notions of Euclidean calculus. In order to
intoduce a "Euclidean derivative", we will take the advantage to have a
distinguished infinite number, namely $\alpha ;$ then we can define a
distinguished infinitesimal number as follows:%
\begin{equation*}
\eta :=\frac{1}{\alpha }.
\end{equation*}

\begin{definition}
\label{md}The \textbf{right derivative} of a normal function $f:\left(
a,b\right) \rightarrow \mathbb{E}$ in a standard point $x_{0}\in (a,b)$, in
the sense of Euclidean Calculus, is defined as follows: 
\begin{equation*}
D^{+}f(x_{0})=ctr\left( \frac{f(x_{0}+\eta )-f(x_{0})}{\eta }\right) ;
\end{equation*}%
similarly the \textbf{left} $\mathbb{E}$\textbf{-derivative} is defined as
follows:%
\begin{equation*}
D^{-}f(x_{0})=ctr\left( \frac{f(x_{0})-f(x_{0}-\eta )}{\eta }\right) ;
\end{equation*}%
the \textbf{mean} $\mathbb{E}$\textbf{-derivative} is defined as follows:%
\begin{equation*}
Df(x_{0})=\frac{1}{2}\left[ D^{+}f(x_{0})+D^{-}f(x_{0})\right] =ctr\left( 
\frac{f(x_{0}+\eta )-f(x_{0}-\eta )}{2\eta }\right) .
\end{equation*}%
We say that a function is derivable in a point $x_{0}\in (a,b)$ if $%
Df(x_{0})=D^{+}f(x_{0})$ and $Df(x_{0})\in \mathbb{R}$. In this case,%
\begin{equation*}
Df(x_{0})=st\left( \frac{f(x_{0}+\eta )-f(x_{0})}{\eta }\right)
\end{equation*}%
is called generalized derivative in the sense of Euclidean Calculus or
simply $\mathbb{E}$\textbf{-derivative.}
\end{definition}

It is easy to check that given a real function $f$ differentiable in a point 
$x_{0}\in \left( a,b\right) \cap \mathbb{R},$ then%
\begin{equation*}
Df(x_{0})=f^{\prime }(x_{0})
\end{equation*}%
but the converse is not true; for example the function $f(x)=x\sin \frac{1}{%
x^{2}}$ is not differentiable for $x=0,$ but 
\begin{equation*}
D^{+}f(0)=ctr\left( \sin \frac{1}{\eta ^{2}}\right) =ctr\left( \sin \frac{1}{%
\left( -\eta \right) ^{2}}\right) =D^{-}f(0)
\end{equation*}%
and hence the $\mathbb{E}$-derivative, given by 
\begin{equation*}
\left[ D\left( x\sin \frac{1}{x^{2}}\right) \right] _{x=0}=st\left( \sin 
\frac{1}{\eta ^{2}}\right) ,
\end{equation*}%
is well defined.

The notion of $\mathbb{E}$-derivative is more general and consequently the
fact that if $f\ $is $\mathbb{E}$-derivable does not imply that $f$
(resticted to $\mathbb{R}$) is continuous. For example the Dirichlet
function 
\begin{equation*}
f_{D}(x):=\left\{ 
\begin{array}{cc}
1 & if\ x\in \mathbb{Q}^{\ast } \\ 
&  \\ 
0 & if\ x\in \mathbb{R}^{\ast }\backslash \mathbb{Q}^{\ast }%
\end{array}%
\right.
\end{equation*}%
is $\mathbb{E}$-derivable in every point with $Df_{D}=0$ (remember that by
Th. \ref{anum}, $\eta \in \mathbb{Q}^{\ast }$), but it is not continuous.

As far, we have defined the derivative of a normal function in a standard
point. The following definition extends the notion of derivative to every
Euclidean point, namely it defines the function "derivative" in all the
points of $\left( a,b\right) $ while Def. \ref{md} defines it only for $x\in
\left( a,b\right) \cap \mathbb{R}$.

\begin{definition}
\label{billa}The $\mathbb{E}$-derivative of a normal derivable function $%
f:\left( a,b\right) \rightarrow \mathbb{E}$ is defined by 
\begin{equation*}
Df:=\left( Df|_{\left( a,b\right) \cap \mathbb{R}}\right) ^{\ast }.
\end{equation*}
\end{definition}

\textbf{Example}: take $f(x)=x\sin \frac{1}{x^{2}};$ then 
\begin{equation*}
D\left( x\sin \frac{1}{x^{2}}\right) =\left\{ 
\begin{array}{cc}
\sin \frac{1}{x^{2}}-\frac{2}{x^{2}}\cos \frac{1}{x^{2}} & if\ \ x\in 
\mathbb{E}\backslash \left\{ 0\right\} \\ 
&  \\ 
st\left( \sin \frac{1}{\eta ^{2}}\right) & if\ \ x=0%
\end{array}%
\right. .
\end{equation*}

Obviously, the $\mathbb{E}$-derivability does not imply the
differentiability defined as follows:

\begin{definition}
A normal function $f:\left( a,b\right) \rightarrow \mathbb{E}$ is said to be
differentiable in a point $x_{0}\in (a,b)$ if there exists a linear function 
$t\mapsto df\left( x_{0}\right) \left[ t\right] $ such that, for every
infinitesimal $\varepsilon $, 
\begin{equation*}
f(x_{0}+\varepsilon )=f(x_{0})+df\left( x_{0}\right) \left[ \varepsilon %
\right] +\varepsilon \varepsilon _{1}
\end{equation*}%
where $\varepsilon _{1}$ is an infinitesimal (which might depend on $%
\varepsilon $).
\end{definition}

It is immediate to check that a function is differentiable in $x_{0}$ if and
and only if%
\begin{equation*}
\forall \varepsilon \in \mathfrak{mon}\left( 0\right) \backslash \left\{
0\right\} ,\ \ Df\left( x_{0}\right) =st\left( \frac{f(x_{0}+\varepsilon
)-f(x_{0})}{\varepsilon }\right)
\end{equation*}

Then, if a function is differentiable in $x_{0}\in (a,b)\cap \mathbb{R}$, it
has the $\mathbb{E}$-derivative in that point and 
\begin{equation*}
df\left( x_{0}\right) \left[ t\right] =Df(x_{0})\cdot t.
\end{equation*}%
but the converse is not true. The derivability of a function does not
coincide with the differentiability; it is well known that in "classic
calculus" this phenomenon occurs only in dimension $\geq 2.$

\begin{remark}
With the above definitions, the two classical problems of the "istantaneous
velocity" and of the "tangent" get different solutions given by the $\mathbb{%
E}$-derivative and the differential respectively. They coincide only for
continuous functions.
\end{remark}

Even if the $\mathbb{E}$-derivative is weaker than the usual one, it is 
\textit{quite surprising} that the main theorems of calculus remain true.
For example, let us consider the Fermat theorem:

\begin{theorem}
(\textbf{Fermat theorem) }If a normal function $f:\left( a,b\right)
\rightarrow \mathbb{E}$ achieves a local maximum (or minimum) in a point $%
x_{0}\in (a+\eta ,b-\eta )$ and it has the $\mathbb{E}$-derivative in that
point, then%
\begin{equation*}
Df(x_{0})=0
\end{equation*}
\end{theorem}

\textbf{Proof}: We have that $D^{+}f(x_{0})\leq 0,\ $and $D^{-}f(x_{0})\geq
0.$ Since $D^{+}f(x_{0})=$ $D^{-}f(x_{0}),$ it follows that%
\begin{equation*}
Df(x_{0})=D^{+}f(x_{0})=0.
\end{equation*}

$\square $

\bigskip

Following the usual procedure, we can prove Rolle theorem, the Lagrange
intermediate value theorem and most of the theorems of real calculus for a
class of function that are not necessarely differentiable, but have only the 
$\mathbb{E}$-derivative. We will sketch this fact (see \cite{bencilibro} for
details). First of all we recall some well known fact in NSA:

\begin{definition}
A function 
\begin{equation*}
f:D\rightarrow \mathbb{E},\ D\subset \mathbb{E}
\end{equation*}%
is called \textbf{continuous in a point} $\xi \in D$ iff 
\begin{equation*}
\xi \sim x\Rightarrow f\left( \xi \right) \sim f\left( x\right) ;
\end{equation*}%
It is called \textbf{continuous in }$D$ if it is normal and it is continuous
in every point $\xi \in D\cap \mathbb{R}$. It is called \textbf{uniformly
continuous in }$D$ if it is normal and it is continuous in every point $\xi
\in D.$
\end{definition}

\bigskip

\begin{theorem}
(\textbf{Weierstrass}) Let $f$ be a continuous fuction on an interval $\left[
a,b\right] .$ Then $f$ has a maximum point in $\left[ a,b\right] \ $and it
is a standard point.
\end{theorem}

\textbf{Proof}: Since the set $\left[ a,b\right] \cap \mathbb{R}%
^{\circledcirc }$ is hyperfinite, $f$ restricted to $\left[ a,b\right] \cap 
\mathbb{R}^{\circledcirc }$ has a maximum point $\xi .$ We calim that $%
c=st\left( \xi \right) $ is the maximum in $\left[ a,b\right] ;$ in fact,
since $\left[ a,b\right] \cap \mathbb{R}\subset \left[ a,b\right] \cap 
\mathbb{R}^{\circledcirc },$ $\forall x\in \left[ a,b\right] \cap \mathbb{R}$%
\begin{equation*}
f\left( \xi \right) \geq f\left( x\right)
\end{equation*}%
and hence by the continuity of $f,$ $\forall x\in \left[ a,b\right] \cap 
\mathbb{R}$%
\begin{equation*}
f\left( c\right) =st\left[ f\left( \xi \right) \right] \geq st\left[ f\left(
x\right) \right] =f(x).
\end{equation*}%
The inequality above, can be extended to every $\zeta =\lim_{\lambda
\uparrow \Lambda }x_{\lambda }\in \left[ a,b\right] ,$ $x_{\lambda }\in %
\left[ a,b\right] \cap \mathbb{R}$. In fact, since $f\left( c\right) \geq
f(x_{\lambda }),$

\begin{equation*}
f\left( c\right) =\lim_{\lambda \uparrow \Lambda }f\left( c\right)
=\lim_{\lambda \uparrow \Lambda }f\left( x_{\lambda }\right) =f\left( \zeta
\right) .
\end{equation*}

$\square $

\bigskip

So we have the following result involving the $\mathbb{E}$-derivative:

\begin{lemma}
(\textbf{Rolle}) Let $f$ be a continuous fuction on an interval $\left[ a,b%
\right] $ such that $f(a)=$ $f(b);$ then if $f$ is $\mathbb{E}$-derivable in 
$\left( a,b\right) $, there is a point $c\in \left( a,b\right) $ such that%
\begin{equation*}
Df\left( c\right) =0
\end{equation*}
\end{lemma}

\textbf{Proof}: By Fermat's and Weierstrass' theorems, the proof is equal to
the usual one.

$\square $

\begin{theorem}
(\textbf{Lagrange}) Let $f$ be a continuous fuction on an interval $\left[
a,b\right] $ and $\mathbb{E}$-derivable in $\left( a,b\right) $, there is a
point $c\in \left( a,b\right) $ such that%
\begin{equation*}
Df\left( c\right) =\frac{f(b)-f(a)}{b-a}.
\end{equation*}
\end{theorem}

\textbf{Proof}: By Rolle's lemma, the proof is equal to the usual one.

$\square $

\bigskip

These results show that even if a $\mathbb{E}$-derivable function can be
quite wild (think of the Dirichlet function), the continuous $\mathbb{E}$%
-derivable functions behave quite well. For example, the space of the
solutions of the equation%
\begin{equation*}
Df=0
\end{equation*}
in general, is not finite-dimensional. However, by the Lagrange theorem it
follows that the only continuous functions which solve the above equations
are the constants. Among the other consequences of the Lagrange's teorem, we
get the following result:

\begin{theorem}
\label{sho}A sufficient condition for a function $f$ to be differentiable in 
$x_{0}\in (a,b)\cap \mathbb{R}$ is that both $f$ and $Df$ be continuous in $%
x_{0}.$
\end{theorem}

\textbf{Proof}: See \cite{bencilibro}.

$\square $

\bigskip

This discussion shows that the notion of $\mathbb{E}$-derivability, even if
it is essentially irrelevant for the applications, it seems interesting for
the foundation of the notion of derivative and its relation with the
differentiability.

\begin{remark}
In the framework of Euclidean calculus there are several other notions of
"generalzed derivative" which make sense. For example we can define the
right \textbf{grid derivative} as follows:%
\begin{equation*}
D^{+}f(x_{0})=ctr\left( \frac{f(x^{+})-f(x)}{x^{+}-x}\right)
\end{equation*}%
where 
\begin{equation*}
x^{+}=\min \ \left\{ y\in \mathbb{R}^{\circledcirc }\ |\ y>x\right\}
\end{equation*}%
and similarly the left grid derivative etc. An other notion of generalized
derivative useful for the applications can be found in \cite{benciU}. In
this paper the notion of grid derivative is combined with the notion of weak
derivative in such a way to include the derivative of distributions
(identified with suitable internal functions).
\end{remark}

\subsection{The integral\label{i}}

Also the definition of the integral takes advantage of a peculiarity of
Euclidean analysis, namely of the operator "$\circledcirc $" introduced by
Def. \ref{tondo}.

\begin{definition}
Given a normal function $f:\left[ a,b\right] \rightarrow \mathbb{E}$, we
define the $\mathbb{E}$-integral as follows:%
\begin{equation*}
\int_{a}^{b}f(x)dx:=ctr\left( \sum_{x\in \left[ a,b\right] ^{\circledcirc
}}f(x)\left( x^{+}-x\right) \right)
\end{equation*}%
where%
\begin{equation*}
x^{+}=\min \ \left\{ y\in \mathbb{R}^{\circledcirc }\ |\ y>x\right\} .
\end{equation*}
\end{definition}

\bigskip

Clearly, if $f$ is Riemann integrable, the $\mathbb{E}$-integral coincides
with the Riemann integral. Moreover the $\mathbb{E}$-integral is well
defined for every normal function even when $\left[ a,b\right] =\mathbb{R}$
or/and $f$ is unbounded. However, the most interesting property of the $%
\mathbb{E}$-integral is given by the following theorem:

\begin{theorem}
Let $f$ be a bounded Lebesgue integrable function, then the $\mathbb{E}$%
-integral is equal to the Lebesgue integral.
\end{theorem}

\textbf{Proof}: Assume that $f$ is a bounded Lebesgue integrable function in 
$\left[ a,b\right] $ and set 
\begin{equation*}
f_{\lambda }(x):=\sum_{z\in \left[ a,b\right] \cap \lambda }f(z)\chi _{\left[
z,z_{\lambda }^{+}\right) }(x)
\end{equation*}%
where $\chi _{\left[ z,z_{\lambda }^{+}\right) }$ is the characteristic
function of $\left[ z,z_{\lambda }^{+}\right) $ and%
\begin{equation*}
z_{\lambda }^{+}=\min \ \left\{ y\in \mathbb{\lambda }\ |\ y>z\right\} .
\end{equation*}%
Now let us denote by $\int_{L}$ the Lebesgue integral; then%
\begin{equation*}
\int_{L}f_{\lambda }(x)dx=\int f_{\lambda }(x)dx
\end{equation*}%
The net of functions $\left\{ f_{\lambda }\right\} $ converges to $f$ in
every point $x\in \left[ a,b\right] \cap \mathbb{R}$ since $\left[ a,b\right]
\cap \mathbb{R}\subset \left[ a,b\right] ^{\circledcirc };$ in fact,
eventually we have that $\forall x\in \left[ a,b\right] \cap \mathbb{R}$ 
\begin{equation*}
f_{\lambda }(x)=f(x)
\end{equation*}%
and hence by (\ref{bn}), $\forall x\in \left[ a,b\right] \cap \mathbb{R}$,%
\begin{equation*}
\lim_{\lambda \rightarrow \Lambda }\ f_{\lambda }(x)=st\left[ f(x)\right]
=f(x)
\end{equation*}%
where $\lim_{\lambda \rightarrow \Lambda }$ is the usual Cauchy limit. By
the Dominated Convergence Theorem, 
\begin{equation*}
\lim_{\lambda \rightarrow \Lambda }\int_{L}f_{\lambda
}(x)dx=\int_{L}\lim_{\lambda \rightarrow \Lambda }\ f_{\lambda
}(x)dx=\int_{L}f(x)dx
\end{equation*}%
On the other hand, 
\begin{eqnarray*}
\lim_{\lambda \uparrow \Lambda }\ \int f_{\lambda }(x)dx &=&\lim_{\lambda
\uparrow \Lambda }\ \left[ \sum_{z\in \left[ a,b\right] \cap \lambda
}f(z)\chi _{\left[ z,z_{\lambda }^{+}\right) }(x)\right] \\
&=&\sum_{x\in \left[ a,b\right] ^{\circledcirc }}f(x)\left( x^{+}-x\right)
\end{eqnarray*}%
Then, by (\ref{bn}), 
\begin{equation*}
\int_{L}f(x)dx\sim \sum_{x\in \left[ a,b\right] ^{\circledcirc }}f(x)\left(
x^{+}-x\right)
\end{equation*}%
and hence 
\begin{equation*}
\int_{L}f(x)dx=st\left( \sum_{x\in \left[ a,b\right] ^{\circledcirc
}}f(x)\left( x^{+}-x\right) \right) =\int f(x)dx.
\end{equation*}

$\square $

\bigskip

\begin{remark}
All the normal functions are $\mathbb{E}$-integrable; however the functions
which are not Lebesgue-integrable might have a pathological behavior; for
example their integral is not invariant for translations. Nevertheless $%
\mathcal{L}^{1}$, the space of the Lebesgue-integrable functions, can be
easily characterized as the closure of the continuous functions with compact
support with respect to the norm%
\begin{equation*}
\left\Vert f\right\Vert =\int \left\vert f(x)\right\vert dx.
\end{equation*}
\end{remark}

\section{Consistency of the axioms}

In this sections, we will prove the consistency of the three axioms
introduced in section \ref{FA} by building a model in ZFC+\{Axiom of
Inaccessibility\}.

\begin{remark}
In many models on nonstandard analysis, the axiom of regularity of ZFC may
fail for external sets (see e.g. \cite{DN,BDNF2}). However, in this paper,
both, the standard universe $V\left( \mathbb{R}\right) $ and the Euclidean
universe $V\left( \mathbb{E}_{\kappa }\right) $ contain only sets of finite
rank and this peculiarity allows to have a model in ZFC. We are forced to
work with sets of finite rank by axiom \ref{num}. In fact, assuming that $%
\Lambda $ contains a set $A$ of infinite rank, we would get a contradiction.
Take for instance a set a defined as follows: 
\begin{equation*}
A=\left\{ b_{n}\ |\ n\in \mathbb{N}\right\}
\end{equation*}%
where $b_{0}=a$ and 
\begin{equation*}
b_{n+1}=\left( b_{n},a\right)
\end{equation*}%
namely 
\begin{equation*}
A=\left\{ a,\left( a,a\right) ,\left( a,a,a\right) ,\left( a,a,a,a\right)
,,....\right\} .
\end{equation*}%
Then, setting 
\begin{equation*}
B=A\times \left\{ a\right\} =\left\{ \left( a,a\right) ,\left( a,a,a\right)
,\left( a,a,a,a\right) ,....\right\} =\left\{ b_{n}\ |\ n\in \mathbb{N}%
^{+}\right\}
\end{equation*}%
we have that $B\subset A$ and hence, by Axiom \ref{num}.\ref{n1} 
\begin{equation*}
\mathfrak{num}\left( B\right) <\mathfrak{num}\left( A\right) ,
\end{equation*}%
while, by Axiom \ref{num}.\ref{n5}%
\begin{equation*}
\mathfrak{num}\left( B\right) =\mathfrak{num}\left( A\times \left\{
a\right\} \right) =\mathfrak{num}\left( A\right) \cdot \mathfrak{num}\left(
\left\{ a\right\} \right) =\mathfrak{num}\left( A\right) \cdot 1=\mathfrak{%
num}\left( A\right) .
\end{equation*}%
Contradiction!
\end{remark}

\subsection{The construction of the field $\mathbb{E}$\label{SE}}

We assume that $\mathbb{A}$ is a set of atoms having cardinality $\left\vert 
\mathbb{A}\right\vert =\kappa $ and that it contains a set $\mathbb{R}$
isomorphic to the real numbers. Moreover, we assume that $\Lambda _{S}$, $%
\Lambda $ and $\mathfrak{L}$ be sets as defined in section \ref{NA}.

If $n_{0}\in \mathbb{N},\ $and $\lambda _{0}\in \mathfrak{L}$, we set 
\begin{equation*}
Q\left( n_{0},\lambda _{0}\right) =\left\{ V_{n}\left( \lambda \right) \in 
\mathfrak{L}\ |\ n\in \mathbb{N},\ n\geq n_{0};\ \lambda \in \mathfrak{L},\
\lambda \supseteq \lambda _{0}\right\}
\end{equation*}%
We have that $Q\left( n_{0},\lambda _{0}\right) \subset \mathfrak{L}$ and%
\begin{equation*}
Q\left( n_{0},\lambda _{0}\right) \cap Q\left( m_{0},\mu _{0}\right)
=Q\left( \max \left\{ n_{0},m_{0}\right\} ,\ \lambda _{0}\cup \mu _{0}\right)
\end{equation*}%
Hence there exists an ultrafilter $\mathcal{U}$ over $\mathfrak{L}$ such
that $\forall n\in \mathbb{N}$, $\forall \lambda \in \mathfrak{L}$, 
\begin{equation}
Q\left( n,\lambda \right) \in \mathcal{U}.  \label{uf}
\end{equation}%
Now, for any $j<\kappa $, we define by transfinite induction a sequence of
ordered fields $\mathbb{K}_{j}\subset \mathbb{A}$ such that $\left\vert 
\mathbb{K}_{j}\right\vert <\kappa $ and 
\begin{equation*}
\mathbb{K}_{j}\subset \mathbb{K}_{k}\ \ \ \ \,\text{if}\ \ j<k.
\end{equation*}%
For $j=0,$ we set 
\begin{equation*}
\mathbb{K}_{0}=\mathbb{R}
\end{equation*}%
and for every ordinal $j<\kappa ,$ we set%
\begin{equation}
\mathfrak{F}_{j+1}:=\mathfrak{F}\left( \mathfrak{L},\mathbb{K}_{j}\right)
/I_{j}  \label{palla}
\end{equation}%
where 
\begin{equation*}
\mathfrak{F}_{j}\left( \mathfrak{L},\mathbb{K}_{j}\right) =\left\{ \varphi
\in \left( \mathbb{K}_{j}\right) ^{\mathfrak{L}}\ |\ \forall \lambda \in 
\mathfrak{L},\ \varphi (\lambda )=\varphi (\lambda \cap \mathbb{K}%
_{j})\right\}
\end{equation*}%
and $I_{j}\subset \mathfrak{F}_{j}\left( \mathfrak{L},\mathbb{K}_{j}\right) $
is the maximal ideal defined as follows:%
\begin{equation*}
I_{j}:=\left\{ \varphi \in \mathfrak{F}_{j}\left( \mathfrak{L},\mathbb{K}%
_{j}\right) \ |\ \exists Q\in \mathcal{U},\ \forall \lambda \in Q,\ \varphi
(\lambda )=0\right\}
\end{equation*}%
Then $\mathfrak{F}_{j+1}$ is a field and the projection%
\begin{equation*}
\Pi _{j}:\mathfrak{F}_{j}\left( \mathfrak{L},\mathbb{K}_{j}\right)
\rightarrow \mathfrak{F}_{j+1}
\end{equation*}%
defined by%
\begin{equation*}
\Pi _{j}\left( \varphi \right) =\left[ \varphi \right] _{\mathcal{U}%
}:=\varphi +I_{j}
\end{equation*}%
is a surjective ring homomorphism. Now, since $\left\vert \mathbb{K}%
_{j}\right\vert <\kappa ,$ also $\left\vert \mathfrak{F}_{j+1}\right\vert
<\kappa ;$ then we can define an injetive map 
\begin{equation*}
\Theta _{j}:\mathfrak{F}_{j+1}\rightarrow \mathbb{A}
\end{equation*}%
such that $\forall \xi \in \mathbb{K}_{j}$, 
\begin{equation*}
\Theta _{j}\left( \left[ C_{\xi }\right] _{\mathcal{U}}\right) =\xi
\end{equation*}%
where $\forall \lambda \in \mathfrak{L,\ }C_{\xi }\left( \lambda \right)
:=\xi \in \mathbb{K}_{j}$.\ The set $\mathbb{K}_{j+1}:=\func{Im}\Theta _{j}$
can be equipped with the structure of ordered field by setting, $\forall \xi
,\zeta \in \mathbb{K}_{j+1}$ 
\begin{eqnarray*}
\xi +\zeta &=&\Theta _{j}\left( \Theta _{j}^{-1}\left( \xi \right) +\Theta
_{j}^{-1}\left( \zeta \right) \right) \\
\xi \cdot \zeta &=&\Theta _{j}\left( \Theta _{j}^{-1}\left( \xi \right)
\cdot \Theta _{j}^{-1}\left( \zeta \right) \right) .
\end{eqnarray*}%
Then $\mathbb{K}_{j+1}$ is a field which contains the real numbers. Now set%
\begin{equation*}
J_{j}=\Theta _{j}\circ \Pi _{j}:\mathfrak{F}\left( \mathfrak{L},\mathbb{K}%
_{j}\right) \rightarrow \mathbb{K}_{j+1}
\end{equation*}%
By this construction, $J_{j}$ is a surjective ring homomorphism.

If $k\leq \kappa $ is a limit ordinal, we set 
\begin{equation*}
\mathbb{K}_{k}=\dbigcup\limits_{j<k}\mathbb{K}_{j}
\end{equation*}%
and 
\begin{equation*}
J_{k}=\ \underset{\longrightarrow }{\lim }J_{j}
\end{equation*}%
namely, $J_{k}\left( \varphi \right) =J_{j}\left( \varphi \right) $ for
every $j$ such that $\varphi \in \mathfrak{F}\left( \mathfrak{L},\mathbb{K}%
_{j}\right) .$ Then, also 
\begin{equation*}
J_{k}:\mathfrak{F}\left( \mathfrak{L},\mathbb{K}_{k}\right) \rightarrow 
\mathbb{K}_{k}
\end{equation*}%
is a surjective ring homomorphism.

In particular, if $j=\kappa $, we have that%
\begin{equation}
J:=J_{\kappa }:\mathfrak{F}\left( \mathfrak{L},\mathbb{K}_{\kappa }\right)
\rightarrow \mathbb{K}_{\kappa }\subset \mathbb{A}  \label{J}
\end{equation}%
is a surjective ring homomorphism and $\left\vert \mathbb{E}_{\kappa
}\right\vert =\kappa $.

If $A\in \Lambda $, by the definition of $\Lambda $ (see (\ref{l1})), we
have that $\left\vert A\right\vert <\kappa ,$ then the net $\left\{ \lambda
\mapsto \left\vert \lambda \cap A\right\vert \right\} \in \mathfrak{F}\left( 
\mathfrak{L},\mathbb{K}_{\kappa }\right) .$ Then, we can define the
numerosity of $A$ as follows:%
\begin{equation}
\mathfrak{num}\left( A\right) =J\left( \left\vert \lambda \cap A\right\vert
\right)  \label{n}
\end{equation}%
where, with some abuse of notation, $\left\vert \lambda \cap A\right\vert $
denotes the net $\left\{ \lambda \mapsto \left\vert \lambda \cap
A\right\vert \right\} .$

\subsection{Proof of the consistency of Axioms 1-3}

Now we can prove the consistency of our axioms.

\begin{theorem}
The numerosity function defined by (\ref{n}) satisfies the request of Axiom %
\ref{num}.
\end{theorem}

\textbf{Proof }: \ref{num}.\ref{n0} and \ref{num}.\ref{n1} follows directly
by the definition (\ref{n}) of $\mathfrak{num}$.

\ref{num}.\ref{n3} - We have that%
\begin{eqnarray*}
\mathfrak{num}\left( A\cup B\right) &=&\ J\left( \left\vert \lambda \cap
\left( A\cup B\right) \right\vert \right) =J\left( \left\vert \left( \lambda
\cap A\right) \cup \left( \lambda \cap B\right) \right\vert \right) \\
&=&J\left( \left\vert \left( \lambda \cap A\right) \right\vert +\left\vert
\left( \lambda \cap B\right) \right\vert \right) =J\left( \left\vert \left(
\lambda \cap A\right) \right\vert \right) +J\left( \left\vert \left( \lambda
\cap A\right) \right\vert \right) \\
&=&\mathfrak{num}\left( A\right) +\mathfrak{num}\left( B\right) .
\end{eqnarray*}

\ref{num}.\ref{n4} - Let $n_{0}$ be so large that $A,B\in V_{n_{0}}\left( 
\mathbb{E}_{\kappa }\right) ,$ then $A\times B\in V_{n_{0}+2}\left( \mathbb{E%
}_{\kappa }\right) ;$ if we take $V_{n}\left( \lambda \right) \in Q\left(
n_{0}+2,\lambda \right) ,$ we have that 
\begin{eqnarray}
\left\vert V_{n}\left( \lambda \right) \cap \left( A\times B\right)
\right\vert &=&\left\vert \left( V_{n}\left( \lambda \right) \cap A\right)
\times \left( V_{n}\left( \lambda \right) \cap B\right) \right\vert  \notag
\\
&=&\left\vert V_{n}\left( \lambda \right) \cap A\right\vert \cdot \left\vert
V_{n}\left( \lambda \right) \cap B\right\vert  \label{abbad}
\end{eqnarray}%
If we set 
\begin{equation*}
\varphi _{E}\left( \lambda \right) =\left\vert V_{n}\left( \lambda \right)
\cap E\right\vert
\end{equation*}%
by (\ref{abbad}), we have that $\forall \lambda \in Q\left( n_{0}+2,\lambda
_{0}\right) ,$%
\begin{equation*}
\varphi _{A\times B}\left( \lambda \right) =\varphi _{A}\left( \lambda
\right) \cdot \varphi _{B}\left( \lambda \right)
\end{equation*}%
and since $Q\left( n_{0}+2,\lambda _{0}\right) \in \mathcal{U}$,%
\begin{equation*}
J\left( \varphi _{A\times B}\right) =J\left( \varphi _{A}\right) \cdot
J\left( \varphi _{B}\right)
\end{equation*}%
Then%
\begin{equation}
\mathfrak{num}\left( A\times B\right) =J\left( \varphi _{A\times B}\right)
=J\left( \varphi _{A}\right) \cdot J\left( \varphi _{B}\right) =\mathfrak{num%
}\left( A\right) \cdot \mathfrak{num}\left( B\right)  \label{5}
\end{equation}

\ref{num}.\ref{n5} - It follows immediately from (\ref{5}).

$\square $

\bigskip

\begin{theorem}
The set $\mathbb{E}_{\kappa }$ satisfies the request of Axiom \ref{1}.
\end{theorem}

\textbf{Proof }: Trivial by our construction.

$\square $

\bigskip

We define the set of the \textbf{Euclidean integers} as follows:%
\begin{equation*}
\mathfrak{Z:=}\dbigcup\limits_{j<\kappa }\mathfrak{Z}_{j}
\end{equation*}%
where 
\begin{equation*}
\mathfrak{Z}_{0}=\mathbb{Z}
\end{equation*}%
and, for $j<\kappa $ 
\begin{equation}
\mathfrak{Z}_{j}:=J\left( \dbigcup\limits_{k<j}\mathfrak{Z}_{k}\right)
\label{jjj}
\end{equation}%
Before proceeding we need the following

\begin{lemma}
\label{m}Every number $\xi \in \mathbb{E}_{\kappa }$ can be decomposed as
follows:%
\begin{equation}
\xi =\zeta +\mu  \label{z}
\end{equation}%
where $\zeta \in \mathfrak{Z}$ and $0\leq \mu <1$.
\end{lemma}

\textbf{Proof}: Let $\xi \in \mathbb{E}_{j}.$ We argue by transfinite
induction over $j.$ If $j=0,$ (\ref{z}) holds with $\zeta \in \mathbb{Z}$.
If $j>0,$ then,%
\begin{equation*}
\xi :=J\left( \left\{ \xi _{\lambda }\right\} \right)
\end{equation*}%
with $\xi _{\lambda }\in \mathbb{E}_{k}$ for some $k<j.$ By our inductive
assumption, 
\begin{equation*}
\xi _{\lambda }=\zeta _{\lambda }+\mu _{\lambda },\ \ \zeta _{\lambda }\in 
\mathfrak{Z\ \ }\text{and\ \ }0\leq \mu _{\lambda }<1.
\end{equation*}%
Then 
\begin{equation*}
\xi :=J\left( \left\{ \zeta _{\lambda }+\mu _{\lambda }\right\} \right)
=J\left( \left\{ \zeta _{\lambda }\right\} \right) +J\left( \left\{ \mu
_{\lambda }\right\} \right)
\end{equation*}%
By (\ref{jjj}), $\zeta :=J\left( \left\{ \zeta _{\lambda }\right\} \right)
\in \mathfrak{Z}$ and if we set $\mu :=J\left( \left\{ \mu _{\lambda
}\right\} \right) ,$ we have that $0\leq \mu <1.\ $Then (\ref{z}) is
satisfied.

$\square $

\bigskip

We now set%
\begin{equation}
\mathfrak{C:=}\left\{ \zeta +r\ |\ \zeta \in \mathfrak{Z},\ r\in \mathbb{R}%
\right\} .  \label{cen}
\end{equation}%
Clearly $\mathfrak{Z}$ and $\mathbb{R}$ are additive subgroups of $\mathbb{E}%
_{\kappa }$ and hence also $\mathfrak{C}$ is an additive group. $\mathfrak{C}
$ is the set of the centers.

\begin{theorem}
The set $\mathfrak{C}$ defined by (\ref{cen}) satisfies the requests of
Axiom \ref{3}.
\end{theorem}

\textbf{Proof.} Since $\mathbb{E}_{\kappa }$ contains the real numbers, if $%
\theta \in \mathbb{E}_{\kappa }$ is bounded, $st(\theta )$ is well defined.
If $\xi \in \mathbb{E}_{\kappa }$, by Lemma \ref{m}, we can write $\xi
=\zeta +\mu ,\ $with $\zeta \in \mathfrak{Z\ }$and\ \ $0\leq \mu <1;$ then
we set 
\begin{equation*}
ctr(\zeta +\mu )=\zeta +st(\mu ).
\end{equation*}%
Then $ctr(\zeta +\mu )\in \mathfrak{C}$.

$\square $

\end{document}